\documentclass[11pt]{article}

\usepackage{caption2}

\usepackage{graphicx}

\usepackage[numbers]{natbib}

\usepackage[all,cmtip]{xy}
\usepackage{setspace}
\usepackage{mathtools, bbm,youngtab}
\usepackage{amscd,mathrsfs}
\usepackage{amssymb,amsmath,dsfont,amsfonts,amsthm}
\usepackage{euscript,enumerate}

\usepackage{stmaryrd}

\usepackage[cbgreek]{textgreek} 

\usepackage{booktabs}
\usepackage{multirow}

\usepackage{hhline} 

\usepackage{cases} 

\usepackage[retainorgcmds]{IEEEtrantools} 

\usepackage{arydshln} 

\usepackage{color}

\DeclareMathAlphabet{\mathpzc}{OT1}{pzc}{m}{it}

\usepackage{epstopdf}

\allowdisplaybreaks 

 \theoremstyle{plain}
 \newtheorem{thm}{Theorem}
 \newtheorem{cor}[thm]{Corollary}

 \theoremstyle{definition}
 
 \newtheorem{rem}[thm]{Remark}
  
 \numberwithin{equation}{section}

\usepackage{hyperref}
\hypersetup{colorlinks=true, linkcolor=blue}
\hypersetup{colorlinks=true,citecolor=blue}

\makeatletter
\DeclareFontFamily{OMX}{MnSymbolE}{}
\DeclareSymbolFont{MnLargeSymbols}{OMX}{MnSymbolE}{m}{n}
\SetSymbolFont{MnLargeSymbols}{bold}{OMX}{MnSymbolE}{b}{n}
\DeclareFontShape{OMX}{MnSymbolE}{m}{n}{
    <-6>  MnSymbolE5
   <6-7>  MnSymbolE6
   <7-8>  MnSymbolE7
   <8-9>  MnSymbolE8
   <9-10> MnSymbolE9
  <10-12> MnSymbolE10
  <12->   MnSymbolE12
}{}
\DeclareFontShape{OMX}{MnSymbolE}{b}{n}{
    <-6>  MnSymbolE-Bold5
   <6-7>  MnSymbolE-Bold6
   <7-8>  MnSymbolE-Bold7
   <8-9>  MnSymbolE-Bold8
   <9-10> MnSymbolE-Bold9
  <10-12> MnSymbolE-Bold10
  <12->   MnSymbolE-Bold12
}{}

\let\llangle\@undefined
\let\rrangle\@undefined
\DeclareMathDelimiter{\llangle}{\mathopen}%
                     {MnLargeSymbols}{'164}{MnLargeSymbols}{'164}
\DeclareMathDelimiter{\rrangle}{\mathclose}%
                     {MnLargeSymbols}{'171}{MnLargeSymbols}{'171}
\makeatother

\title{\Large{\textbf{A Conformal Three-Field Formulation for Nonlinear Elasticity:\\ From Differential Complexes to Mixed Finite Element Methods}}}

\author{Arzhang Angoshtari\thanks{Department of Civil and Environmental Engineering, The George Washington University, Washington, DC 20052. E-mail: aangoshtari@gwu.edu.} \and 
 Ali Gerami Matin\thanks{Department of Civil and Environmental Engineering, The George Washington University, Washington, DC 20052. E-mail: agerami@gwu.edu.}}

\begin{document}

\maketitle

\begin{abstract} We introduce a new class of mixed finite element methods for $2$D and $3$D compressible nonlinear elasticity. The independent unknowns of these conformal methods are displacement, displacement gradient, and the first Piola-Kirchhoff stress tensor. The so-called edge finite elements of the $\mathrm{curl}$ operator is employed to discretize the trial space of displacement gradients. Motivated by the differential complex of nonlinear elasticity, this choice guarantees that discrete displacement gradients satisfy the Hadamard jump condition for the strain compatibility. We study the stability of the proposed mixed finite element methods by deriving some inf-sup conditions. By considering $32$ choices of simplicial conformal finite elements of degrees $1$ and $2$, we show that $10$ choices are not stable as they do not satisfy the inf-sup conditions. We numerically study the stable choices and conclude that they can achieve optimal convergence rates. By solving several $2$D and $3$D numerical examples, we show that the proposed methods are capable of providing accurate approximations of strain and stress.      
\end{abstract}

\begin{description}
\item[Keywords.] Nonlinear elasticity; mixed finite element methods; inf-sup conditions; differential complex; finite element exterior calculus.
\end{description}


\section{Introduction}

Although modeling deformations of nonlinearly elastic solids is an old problem \citep{TrNo65}, designing stable computational methods for predicting nonlinear deformations in some modern engineering applications such as electroactive polymers and biological tissues is still a challenging task. A simple strategy is to extend well-performing computational methods of linearized elasticity to nonlinear elasticity, however, it is well-known that due to the occurrence of various unphysical instabilities, such extensions may have a very poor performance \citep{AuDaLoReTaWr2013,AuDaLoRe2005}. 

It was shown that mixed finite element methods provide a useful framework for studying compressible and incompressible nonlinear elasticity \citep{Wr2009}. Mixed finite element methods involve several independent unknowns and are usually defined as finite element methods which are based on a primal-dual problem or a saddle-point variational problem \citep{BoBrFo2013}. A different definition of mixed methods is methods that simultaneously approximate an unknown and some of its derivatives \citep{Od1972}, see also \citep[Chapter 7]{Ciarlet1978} for a general classification of finite element methods.           

Different mixed methods exist for compressible nonlinear elasticity; For example, two-field methods based on the Hellinger-Reissner principle in terms of displacement and stress \citep{Wr2009} and three-field methods based on the Hu-Washizu principle in terms of displacement, strain, and stress \citep{SiAr1992,AnFSYa2017}. Some potential advantages of mixed finite element methods in nonlinear elasticity include locking free behavior for thin solids and for the (near) incompressible regime, good performance for problems involving large bending, accurate approximations of strains and stresses, insensitivity to mesh distortions, and simple implementation of constitutive equations \citep{Wr2009,Re2015}. On the other hand, a disadvantage of mixed methods is that they are computationally more expensive comparing to standard single-field methods as there are more degrees of freedom per element. Another disadvantage of mixed methods is that the well-posedness of the underlying mixed formulation is not necessarily inherited by its discretizations. This aspect of mixed methods is usually studied in the context of inf-sup conditions \citep{BoBrFo2013}.

In this paper, we introduce a new class of conformal mixed finite element methods for $2$D and $3$D compressible nonlinear elasticity based on a three-field formulation in terms of displacement, displacement gradient, and the first Piola-Kirchhoff stress tensor. The main idea is to discretize the trial space of displacement gradients by employing finite elements suitable for the $\mathrm{curl}$ operator. This choice can be readily justified by using a mathematical structure called the differential complex of nonlinear elasticity \citep{AngoshtariYavari2014I,AngoshtariYavari2014II} and guarantees that the Hadamard jump condition for the strain compatibility is satisfied on the discrete level as well. A relation between the nonlinear elasticity complex and a well-known complex from differential geometry called the de Rham complex allows one to discretize the former by using finite element spaces suitable for the discretization of the latter. We employ these finite element spaces to derive finite element methods for nonlinear elasticity. We show that even for hyperelastic materials, the underlying weak form does not correspond to a saddle-point problem. However, the resulting finite element methods are still called mixed methods in the sense that displacement and its derivative are approximated simultaneously.

Similar ideas were employed in \citep{AnFSYa2017} to obtain a class of mixed finite element methods for $2$D compressible nonlinear elasticity. In contrary to the present work, one can show that the underlying weak form of \citep{AnFSYa2017} is associated to a saddle-point of a Hu-Washizu-type functional for hyperelastic matrerials. Numerical examples suggested that the resulting finite element methods have good features such as optimal convergence rates, good bending performance, accurate approximations of strains and stresses, and the lack of the hourglass instability that may occur in non-conformal enhanced strain mixed methods \citep{WrRe1996}. However, those mixed methods suffer from at least two drawbacks: On the one hand, only a limited number of finite element choices lead to a stable method and on the other hand, and more importantly, their extension to $3$D problems is hard.

In comparison to \citep{AnFSYa2017}, the present mixed methods work well for both $2$D and $3$D problems and also are stable for broader choices of elements. For example, in \citep{AnFSYa2017} it was observed that only $7$ out of the $32$ possible choices of first-order and second-order triangular elements lead to stable methods while $22$ choices are stable in the present work. The main difference is among the test spaces of the constitutive relation: While curl-based spaces are used in the previous work, divergence-based spaces are employed in this work. It is not hard to see that in the formulation of \citep{AnFSYa2017}, instead of seeking stresses in divergence-based spaces, they are implicitly sought in the intersection of curl-based and divergence-based spaces. The present formulation does not impose this unphysical restriction on stress.

We employ the general framework of \citep{PoRa1994,CaRa1997} for the Galerkin approximation of regular solutions of nonlinear problems to study the stability of the proposed methods. In particular, we write a sufficient inf-sup condition and two other weaker inf-sup conditions. By considering $32$ choices of simplicial finite elements of degrees $1$ and $2$ in $2$D and $3$D, the performance of mixed methods are studied. We show that $10$ choices are not stable as they violate the inf-sup conditions. Our numerical examples suggest that the proposed mixed methods are capable of attaining optimal convergence rates and approximate strains and stresses accurately.

This paper is organized as follows: In Section \ref{Sec_MixFEM}, we first briefly review the differential complex of nonlinear elasticity and then we introduce a mixed formulation for nonlinear elasticity. This formulation is then discretized by employing suitable conformal finite element spaces. In Section \ref{Sec_Conv}, a convergence analysis for regular solutions is presented and suitable inf-sup conditions are written.  Also we rigorously show that some choices of simplicial finite elements do not satisfy these inf-sup conditions. By considering several $2$D and $3$D numerical examples, the performance of the proposed finite element methods is studied in Section \ref{Sec_Examples}. We present a numerical study of the inf-sup conditions as well. Some final remarks will be made in Section \ref{Sec_Conc}.

\section{A Class of Mixed Finite Element Methods for Nonlinear Elasticity}\label{Sec_MixFEM}

The mixed finite element methods introduced in this work are closely related to the nonlinear elasticity complex. We begin with a brief description of this complex for $3$D nonlinear elasticity and then, we employ this complex to introduce a mixed formulation for nonlinear elasticity. Mixed finite element methods are then defined by discretizing this mixed formulation by using suitable conformal finite element spaces. We assume $\{\mathrm{X}^{I}\}_{I=1}^{n}$ and $\{\mathbf{E}^{I}\}_{I=1}^{n}$ are respectively the Cartesian coordinates and the standard orthonormal basis of $\mathbb{R}^{n}$, $n=2,3$. Since covariant and contravariant components of tensors are the same in $\{\mathrm{X}^{I}\}_{I=1}^{n}$, we will only use contravariant components of tensors. Also unless stated otherwise, we use the summation convention on repeated indices. 

\subsection{The Nonlinear Elasticity Complex}
Let $\mathcal{B}$ represent the reference configuration of a $3$D elastic body with the boundary $\partial \mathcal{B}$. The unit outward normal vector field of $\partial\mathcal{B}$ is denoted by $\boldsymbol{N}$ and we assume $\partial\mathcal{B}=\Gamma_{1}\cup\Gamma_{2}$, where $\Gamma_{1}$ and $\Gamma_{2}$ have disjoint interiors. A second-order tensor field $\boldsymbol{T}$ on $\mathcal{B}$ is said to be normal to $\Gamma_{i}$, $i=1,2$, if $\boldsymbol{T}(\mathbf{Y}):=T^{IJ}Y^{J}\mathbf{E}_{I}=\boldsymbol{0}$, for any vector $\mathbf{Y}$ parallel to $\Gamma_{i}$. Similarly, $\boldsymbol{T}$ is said to be parallel to $\Gamma_{i}$ if $\boldsymbol{T}(\mathbf{Y})=\boldsymbol{0}$, for any vector $\mathbf{Y}$ normal to $\Gamma_{i}$.

Given a vector field $\boldsymbol{U}$ and a tensor field $\boldsymbol{T}$, one can define the operators $\mathbf{grad}$, $\mathbf{curl}$, and $\mathbf{div}$ as 
\begin{equation*}
(\mathbf{grad}\,\boldsymbol{U})^{IJ}=\partial_{J}U^{I},~~ (\mathbf{curl}\,\boldsymbol{T})^{IJ}=\varepsilon_{JKL}\partial_{K}T^{IL}, ~~(\mathbf{div}\,\boldsymbol{T})^{I}=\partial_{J}T^{IJ},
\end{equation*}
where ``$\partial_{J}$'' denotes $\partial/\partial \mathrm{X}^{J}$ and $\varepsilon_{JKL}$ is the standard permutation symbol. Suppose $[H^{1}(\mathcal{B})]^{3}$ is the standard space of $H^{1}$ vector fields on $\mathcal{B}$ (i.e. the space of vector fields such that their components and first derivatives of their components are square integrable) and let $[H^{1}_{i}(\mathcal{B})]^{3}$ be the space of $H^{1}$ vector fields that vanish on $\Gamma_{i}$. By $H^{\mathbf{c}}(\mathcal{B})$, we denote the space of second-order tensor fields $\boldsymbol{T}$ such that both $\boldsymbol{T}$ and $\mathbf{curl}\, \boldsymbol{T}$ are of $L^{2}$-class (i.e. have square integrable components). The space of $H^{\mathbf{c}}$ tensor fields that are normal to $\Gamma_{i}$ is denoted by $H^{\mathbf{c}}_{i}(\mathcal{B})$. Similarly, the space of $L^{2}$ second-order tensor fields with $L^{2}$ divergence is denoted by $H^{\mathbf{d}}(\mathcal{B})$ and $H^{\mathbf{d}}_{i}(\mathcal{B})$ indicates $H^{\mathbf{d}}$ tensor fields that are parallel to $\Gamma_{i}$.

It is possible to define continuous operators $\mathbf{grad}: [H^{1}_{i}(\mathcal{B})]^{3}\rightarrow H^{\mathbf{c}}_{i}(\mathcal{B})$, $\mathbf{curl}: H^{\mathbf{c}}_{i}(\mathcal{B})\rightarrow H^{\mathbf{d}}_{i}(\mathcal{B})$, and $\mathbf{div}: H^{\mathbf{d}}_{i}(\mathcal{B}) \rightarrow [L^{2}(\mathcal{B})]^{3}$ that satisfy the relations $\mathbf{curl}(\mathbf{grad}\,\boldsymbol{Y})=\boldsymbol{0}$, and $\mathbf{div}(\mathbf{curl}\,\boldsymbol{T})=\boldsymbol{0}$. These facts are usually expressed by writing the differential complex
\begin{equation}\label{gcdTenHilb3D}
\scalebox{1}{\xymatrix@C=3.5ex{0 \ar[r] & [H^{1}_{i}(\mathcal{B})]^{3} \ar[r]^-{\mathbf{grad}} &H^{\mathbf{c}}_{i}(\mathcal{B}) \ar[r]^-{\mathbf{curl}} & H^{\mathbf{d}}_{i}(\mathcal{B}) \ar[r]^-{\mathbf{div}}  &[L^{2}(\mathcal{B})]^{3} \ar[r]  &0.  } } 
\end{equation}
The above complex is called the nonlinear elasticity complex as it describes the kinematics and the kinetics of nonlinearly elastic bodies in the following sense \citep{AngoshtariYavari2014I, AngoshtariYavari2014II}: Let $\varphi:\mathcal{B}\rightarrow\mathbb{R}^{3}$ be a deformation of $\mathcal{B}$ and let $\boldsymbol{U}$ be the associated displacement field. Then, the displacement gradient is $\boldsymbol{K}:=\mathbf{grad}\,\boldsymbol{U}$, and $\mathbf{curl}\,\boldsymbol{K}=\boldsymbol{0}$, is the necessary condition for the compatibility of $\boldsymbol{K}$. On the other hand, $\mathbf{div}\,\boldsymbol{P}=\boldsymbol{0}$, is the equilibrium equation in terms of the first Piola-Kirchhoff stress tensor $\boldsymbol{P}$. This equation is also the necessary condition for the existence of a stress function $\boldsymbol{\Psi}$ such that $\boldsymbol{P}=\mathbf{curl}\,\boldsymbol{\Psi}$. By considering the $2$D curl operator $(\mathbf{curl}\,\boldsymbol{T})^{I}=\partial_{1}T^{I2}-\partial_{2}T^{I1}$, one can also write similar results for $2$D nonlinear elasticity. One can show that \eqref{gcdTenHilb3D} provides a connection between solutions of certain partial differential equations and the topologies of $\mathcal{B}$ and $\Gamma_{i}$.

\subsection{A Three-Field Mixed Formulation}

Motivated by the complex \eqref{gcdTenHilb3D}, we write a mixed formulation for nonlinear elasticity in terms of the displacement $\boldsymbol{U}$, the displacement gradient $\boldsymbol{K}$, and the first Piola-Kirchhoff stress tensor $\boldsymbol{P}$. Let $\boldsymbol{P}=\mathbb{P}(\boldsymbol{K})$ express the constitutive equation of the elastic body $\mathcal{B}\subset\mathbb{R}^{n}$, $n=2,3$. The boundary value problem of nonlinear elastostatics can be written as: Given a body force $\boldsymbol{B}$, a displacement $\overline{\boldsymbol{U}}$ of $\mathcal{B}$, and a traction vector field $\overline{\boldsymbol{T}}$ on $\Gamma_{2}$, find $(\boldsymbol{U},\boldsymbol{K},\boldsymbol{P})$ such that

\begin{subequations}\label{NonLinStrong}
\begin{IEEEeqnarray}{rlrll}
& \mathbf{div}\,\boldsymbol{P} = - \boldsymbol{B}, & &\quad & \label{NonLinStrong1}
\\*
&\boldsymbol{K}-\mathbf{grad}\,\boldsymbol{U}=\boldsymbol{0},  &\smash{\left.
\IEEEstrut[8\jot]
\right\}} & & \text{in } \mathcal{B}, \label{NonLinStrong2}
\\*
& \boldsymbol{P}-\mathbb{P}(\boldsymbol{K})=\boldsymbol{0}, & & & \label{NonLinStrong3}
\\*
& \boldsymbol{U}=\overline{\boldsymbol{U}}, & & &\text{on } \Gamma_{1}, \label{NonLinStrong4}
\\*
& \boldsymbol{P}(\boldsymbol{N})=\overline{\boldsymbol{T}}, & & &\text{on } \Gamma_{2}. \label{NonLinStrong5}
\end{IEEEeqnarray}
\end{subequations}

To write a weak formulation for the above problem, we proceed as follows: Let ``$\boldsymbol{\cdot}$'' denote the standard inner product of $\mathbb{R}^{n}$ and let $\llangle,\rrangle$ denote both the $L^{2}$-inner product of vector fields $\llangle\boldsymbol{Y},\boldsymbol{Z}\rrangle:=\int_{\mathcal{B}}Y^{I}Z^{I}dV$, and the $L^{2}$-inner product of tensor fields $\llangle\boldsymbol{S},\boldsymbol{T}\rrangle:=\int_{\mathcal{B}}S^{IJ}T^{IJ}dV$. By taking the $L^{2}$-inner product of \eqref{NonLinStrong1} with an arbitrary $\boldsymbol{\Upsilon}\in [H^{1}_{1}(\mathcal{B})]^{n}$ and using Green's formula
\begin{equation*}  
\llangle \boldsymbol{P},\mathbf{grad}\,\boldsymbol{\Upsilon} \rrangle = - \llangle \mathbf{div}\,\boldsymbol{P},\boldsymbol{\Upsilon} \rrangle + \int_{\partial\mathcal{B}} \boldsymbol{P}(\boldsymbol{N})\boldsymbol{\cdot} \boldsymbol{\Upsilon} dA,
\end{equation*}
one concludes that 
\begin{equation*}  
\llangle \boldsymbol{P},\mathbf{grad}\,\boldsymbol{\Upsilon} \rrangle = \llangle \boldsymbol{B},\boldsymbol{\Upsilon} \rrangle + \int_{\Gamma_{2}} \boldsymbol{P}(\boldsymbol{N})\boldsymbol{\cdot} \boldsymbol{\Upsilon} dA. 
\end{equation*}
We also take the $L^{2}$-inner product of \eqref{NonLinStrong2} and \eqref{NonLinStrong3} with arbitrary $\boldsymbol{\lambda}$ of $H^{\mathbf{c}}$-class and arbitrary $\boldsymbol{\pi}$ of $H^{\mathbf{d}}$ class and obtain the following mixed formulation for nonlinear elastostatics:

\bigskip
\begin{minipage}{.95\textwidth}{\textit{Given a body force $\boldsymbol{B}$, a displacement $\overline{\boldsymbol{U}}$ of $\mathcal{B}$, and a boundary traction vector field $\overline{\boldsymbol{T}}$ on $\Gamma_{2}$, find $(\boldsymbol{U},\boldsymbol{K},\boldsymbol{P})\in [H^{1}(\mathcal{B})]^{n}\times H^{\mathbf{c}}(\mathcal{B}) \times H^{\mathbf{d}}(\mathcal{B})$ such that $\boldsymbol{U}= \overline{\boldsymbol{U}}$, on $\Gamma_1$ and }}
\begin{equation}\label{3DElasMixF}
\begin{alignedat}{3}
\llangle\boldsymbol{P},\mathbf{grad}\,\boldsymbol{\Upsilon}\rrangle&=\llangle \boldsymbol{B},\boldsymbol{\Upsilon} \rrangle + \int_{\Gamma_{2}}\overline{\boldsymbol{T}}\boldsymbol{\cdot}\boldsymbol{\Upsilon} dA, &\quad &\forall \boldsymbol{\Upsilon}\in [H^{1}_{1}(\mathcal{B})]^{n},\\
\llangle\mathbf{grad}\,\boldsymbol{U},\boldsymbol{\lambda}\rrangle - \llangle \boldsymbol{K},\boldsymbol{\lambda}\rrangle &= 0,& &\forall \boldsymbol{\lambda}\in H^{\mathbf{c}}(\mathcal{B}),\\
\llangle \mathbb{P}(\boldsymbol{K}),\boldsymbol{\pi}\rrangle - \llangle \boldsymbol{P}, \boldsymbol{\pi} \rrangle  &= 0, & &\forall \boldsymbol{\pi}\in H^{\mathbf{d}}(\mathcal{B}).
\end{alignedat}
\end{equation}
\end{minipage}
\bigskip

\begin{rem} The mixed formulation \eqref{3DElasMixF} is different from that of \citep[Equation (2.8)]{AnFSYa2017}: Here, test functions associated to the definition of the displacement gradient and the constitutive relation are respectively of classes $H^{\mathbf{c}}$ and $H^{\mathbf{d}}$. In \citep{AnFSYa2017}, $H^\mathbf{c}$ test functions are employed for the constitutive relation and $H^{\mathbf{d}}$ test functions for the definition of the displacement gradient. Later we will show that the mixed formulation of \citep{AnFSYa2017} imposes an unphysical constraint on stresses. 
\end{rem}

\begin{rem} For hyperelastic materials, the mixed formulation of \citep{AnFSYa2017} is a saddle-point problem associated to a Hu-Washizu-type functional. However, the mixed formulation \eqref{3DElasMixF} does \emph{not} correspond to a stationary point of any functional $J:Z\rightarrow \mathbb{R}$ with $Z=[H^{1}(\mathcal{B})]^{n}\times H^{\mathbf{c}}(\mathcal{B}) \times H^{\mathbf{d}}(\mathcal{B})$. To show this, let $u,v,w\in Z$, where $u=(\boldsymbol{U},\boldsymbol{K},\boldsymbol{P})$, $v=(\boldsymbol{\Upsilon},\boldsymbol{\lambda},\boldsymbol{\pi})$, $w=(\boldsymbol{V},\boldsymbol{M},\boldsymbol{Q})$, and notice that the problem \eqref{3DElasMixF} can be written as: Find $u\in Z$ such that $G(u,v)=0$, $\forall v\in Z$, where
\begin{equation*}
\begin{alignedat}{3}
G(u,v) &= \llangle\boldsymbol{P},\mathbf{grad}\,\boldsymbol{\Upsilon}\rrangle + \llangle\mathbf{grad}\,\boldsymbol{U},\boldsymbol{\lambda}\rrangle - \llangle \boldsymbol{K},\boldsymbol{\lambda}\rrangle\\
 &+ \llangle \mathbb{P}(\boldsymbol{K}),\boldsymbol{\pi}\rrangle - \llangle \boldsymbol{P}, \boldsymbol{\pi} \rrangle  - \llangle \boldsymbol{B},\boldsymbol{\Upsilon} \rrangle - \int_{\Gamma_{2}}\overline{\boldsymbol{T}}\boldsymbol{\cdot}\boldsymbol{\Upsilon} dA. 
\end{alignedat}
\end{equation*}
If \eqref{3DElasMixF} corresponds to a stationary point of $J:Z\rightarrow \mathbb{R}$, then $G(u,v)=\mathrm{D}J(u)v$, where $\mathrm{D}J(u)v$ is the (Fr\'{e}chet) derivative of $J$ at $u$ in the direction of $v$. Since the second derivative of $J$ has the symmetry $\mathrm{D}^{2}J(u)(v,w)=\mathrm{D}^{2}J(u)(w,v)$ \citep{Av1986}, the first derivative of $G$ should satisfy $\mathrm{D}_{1}G(u,v)w=\mathrm{D}_{1}G(u,w)v$, with 
\begin{equation*}
\mathrm{D}_{1}G(u,v)w = \llangle\boldsymbol{Q},\mathbf{grad}\,\boldsymbol{\Upsilon}\rrangle + \llangle\mathbf{grad}\,\boldsymbol{V},\boldsymbol{\lambda}\rrangle - \llangle \boldsymbol{M},\boldsymbol{\lambda}\rrangle + \llangle \mathsf{A}(\boldsymbol{K})\boldsymbol{:}\boldsymbol{M},\boldsymbol{\pi}\rrangle - \llangle \boldsymbol{Q}, \boldsymbol{\pi} \rrangle,
\end{equation*}
where $\mathsf{A}(\boldsymbol{K})$ is the elasticity tensor in terms of the displacement gradient and $(\mathsf{A}(\boldsymbol{K})\boldsymbol{:}\boldsymbol{M})^{IJ}:=A^{IJRS}M^{RS}$. However, it is easy to check that $\mathrm{D}_{1}G(u,v)w\neq\mathrm{D}_{1}G(u,w)v$, and therefore, the formulation \eqref{3DElasMixF} does not correspond to a saddle-point of any functional. Despite this fact, we still call a finite element method based on \eqref{3DElasMixF} a mixed method in the general sense that displacement and its derivative are independent unknowns of this formulation; See the discussion of \citep[Page 417]{Ciarlet1978} regarding definitions of mixed methods.    
\end{rem}

\begin{rem}\label{NotWellDefined} In general, the response function $\mathbb{P}$ of a nonlinearly elastic body $\mathcal{B}$ is not a well-defined mapping $\mathbb{P}:H^{\mathbf{c}}(\mathcal{B})\rightarrow H^{\mathbf{d}}(\mathcal{B})$, i.e. $\boldsymbol{K}\in H^{\mathbf{c}}(\mathcal{B})$ does not necessarily imply that $\mathbb{P}(\boldsymbol{K})\in H^{\mathbf{d}}(\mathcal{B})$ \citep{An2018}. Roughly speaking, this means that it is impossible to induce arbitrary continuous deformations in nonlinearly elastic bodies by using external loads. To study the well-posedness of the problem \eqref{3DElasMixF} by using approaches based on the implicit function theorem, it is sufficient to assume that the restriction $\mathbb{P}:O\rightarrow [L^{2}(\mathcal{B})]^{3}$ is a well-defined mapping, where $O$ is an open subset of $H^{\mathbf{c}}(\mathcal{B})$.    
\end{rem}

\subsection{Mixed Finite Element Methods}\label{Sec_CSFEM}
The complex \eqref{gcdTenHilb3D} has a close relation with a well-known complex from differential geometry called the de Rham complex \citep{AngoshtariYavari2014II}. One can employ this relation to obtain conformal mixed finite element methods for approximating solutions of \eqref{3DElasMixF} as follows. It was shown that the finite element exterior calculus (FEEC) provides a systematic method for discretizing the de Rham complex using finite element spaces \citep{Arnold2010,Arnold2006}. The relation between \eqref{gcdTenHilb3D} and the de Rham complex allows one to obtain $H^{\mathbf{c}}$- and $H^{\mathbf{d}}$-conformal finite element spaces by using FEEC. For example, to obtain $H^{\mathbf{c}}$-conformal finite element spaces over a $3$D body $\mathcal{B}$, we proceed as follows: Let the (row) vector field $\boldsymbol{K}_{I}=(K^{I1},K^{I2},K^{I3})$ denote the $I$-th row of the displacement gradient $\boldsymbol{K}$. One can write
\begin{equation*}
\mathbf{curl}\,\boldsymbol{K}=\mathbf{curl}\left[ \begin{array}{c} \boldsymbol{K}_{1} \\ \boldsymbol{K}_{2} \\ \boldsymbol{K}_{3} \end{array}\right] = \left[ \begin{array}{c} \mathrm{curl}\,\boldsymbol{K}_{1} \\ \mathrm{curl}\,\boldsymbol{K}_{2} \\ \mathrm{curl}\,\boldsymbol{K}_{3} \end{array}\right],
\end{equation*}         
where $\mathrm{curl}$ is the standard curl operator of vector fields. Consequently, $H^{\mathbf{c}}(\mathcal{B})$ can be identified with three copies of the standard curl space $H^{c}(\mathcal{B})$ for vector fields, i.e. $H^{\mathbf{c}}(\mathcal{B})=[H^{c}(\mathcal{B})]^{3}$. On the other hand, the space $H^{c}(\mathcal{B})$ of vector fields can be identified with a space of differential $1$-forms, which can be discretized using FEEC. Thus, conformal finite element spaces for $H^{\mathbf{c}}(\mathcal{B})$ can be obtained by using three copies of conformal finite element spaces of differential $1$-forms. Similarly, since $H^{\mathbf{d}}(\mathcal{B})=[H^{d}(\mathcal{B})]^{3}$, where $H^{d}(\mathcal{B})$ is the divergence space of vector fields, one can obtain conformal finite element spaces for $H^{\mathbf{d}}(\mathcal{B})$ by using three copies of conformal finite element spaces for $2$-forms. In $2$D, by noting that $H^{\mathbf{c}}(\mathcal{B})=[H^{c}(\mathcal{B})]^{2}$ and $H^{\mathbf{d}}(\mathcal{B})=[H^{d}(\mathcal{B})]^{2}$, one can obtain conformal finite element spaces for $H^{\mathbf{c}}(\mathcal{B})$ and $H^{\mathbf{d}}(\mathcal{B})$ by using two copies of conformal finite element spaces for $1$-forms.

\begin{figure}[t]
\begin{center}
\includegraphics[scale=.7,angle=0]{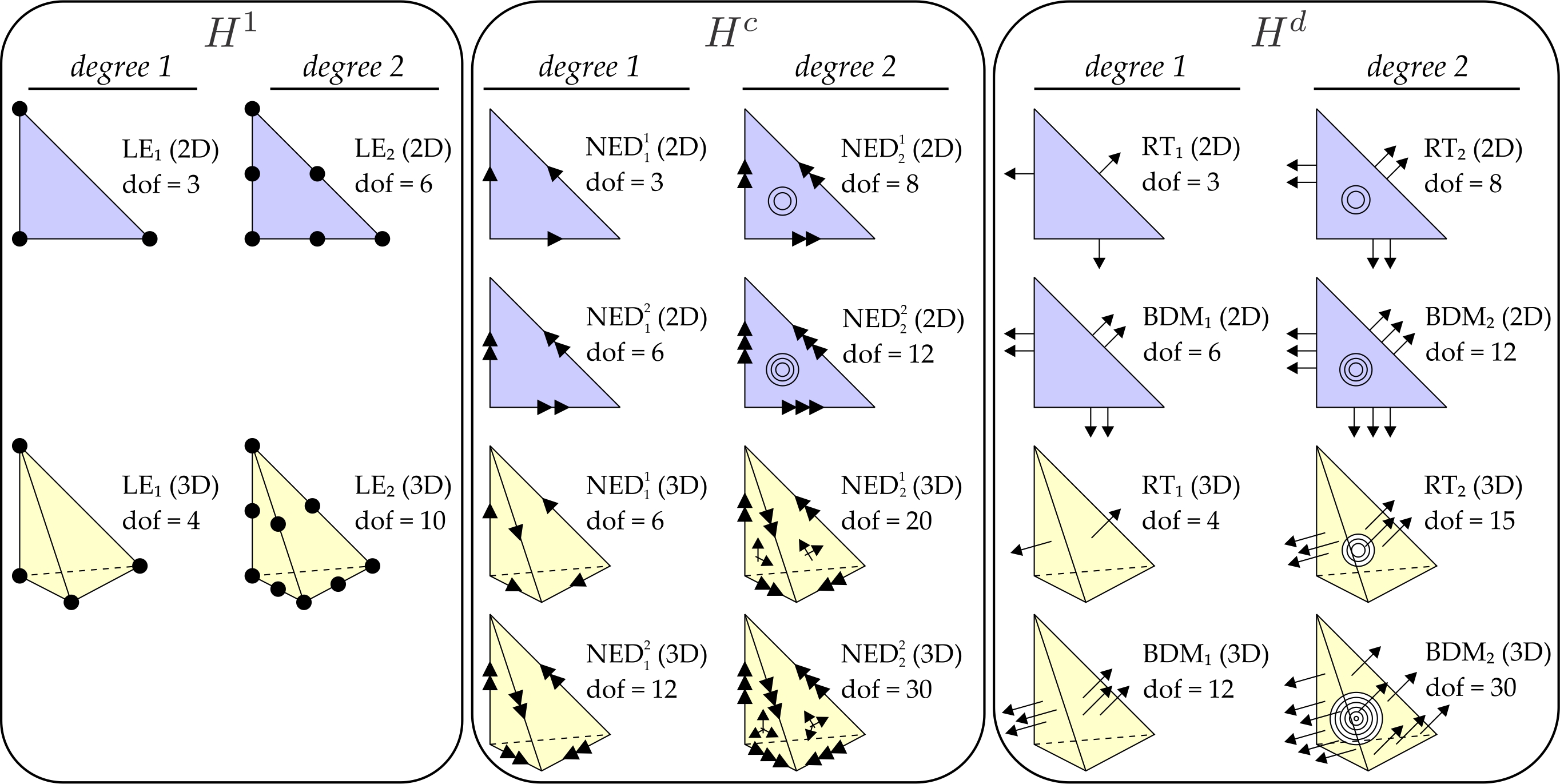}
\end{center}
\vspace*{-0.2in}
\caption{\footnotesize Conventional finite element diagrams of the first and the second degree $H^{1}$, $H^{c}$, and $H^{d}$ elements on triangles and tetrahedra. In this figure, $\mathrm{LE}_i$ stands for the Lagrange element of degree $i$, $\mathrm{NED}^{j}_i$ stands for the $i$-th degree N\'{e}d\'{e}lec element of the $j$-th kind, $\mathrm{RT}_i$ stands for the Raviart-Thomas element of degree $i$, and $\mathrm{BDM}_i$ stands for the Brezzi-Douglas-Marini element of degree $i$. Arrows parallel (normal) to an edge or a face denote degrees of freedom associated to tangent (normal) components of vector fields along that edge or face. Only degrees of freedom associated to visible edges and faces are shown.} 
\label{FEO12}
\end{figure}

Let $\mathcal{B}$ be a polyhedral domain with a simplicial mesh $\mathcal{B}_{h}$, i.e. $\mathcal{B}_{h}$ is a triangular mesh in $2$D and a tetrahedral mesh in $3$D. The above discussion implies that one can associate a tensorial finite element space $V^{\mathbf{c}}_{h}:=[V^{c}_{h}]^{n}\subset H^{\mathbf{c}}(\mathcal{B})$ to any vectorial finite element space $V^{c}_{h}\subset H^{c}(\mathcal{B})$ with $\dim V^{\mathbf{c}}_{h}=n\dim V^{c}_{h}$. Similarly, for any $H^{d}$-conformal finite element space $V^{d}_{h}\subset H^{d}(\mathcal{B})$, one obtains $H^{\mathbf{d}}$-conformal finite element space $V^{\mathbf{d}}_{h}:=[V^{d}_{h}]^{n}\subset H^{\mathbf{d}}(\mathcal{B})$ with $\dim V^{\mathbf{d}}_{h}=n\dim V^{d}_{h}$. FEEC provides a systematic approach for obtaining finite element spaces $V^{c}_{h}$ and $V^{d}_{h}$ of arbitrary order. For example, Figure \ref{FEO12} shows conventional finite element diagrams of some $H^{1}$-, $H^{c}$-, and $H^{d}$-conformal elements of degrees 1 and 2. Notice that for $H^{c}$ elements, some degrees of freedom are associated to tangent components of vectors fields along faces and edges, whereas degrees of freedom of $H^{d}$ elements are associated to normal components of vector fields along faces; See \citep[Chapter 3]{loMaWe2012} for more details about these elements.             

Let $[V^{1}_{h}]^{n}$, $V^{\mathbf{c}}_{h}$, $V^{\mathbf{d}}_{h}$ be finite element spaces as described above and let $V^{1}_{h,i}=V^{1}_{h}\cap H^{1}_{i}(\mathcal{B})$. Also suppose $\mathcal{I}^{1}_{h}$ is the canonical interpolation operators associated to the $H^{1}$ elements. Then, we consider the following mixed finite element methods for \eqref{3DElasMixF}:

\bigskip
\begin{minipage}{.95\textwidth}{\textit{Given a body force $\boldsymbol{B}$, a displacement $\overline{\boldsymbol{U}}$, and a boundary traction vector field $\overline{\boldsymbol{T}}$ on $\Gamma_{2}$, find $(\boldsymbol{U}_{h},\boldsymbol{K}_{h},\boldsymbol{P}_{h})\in [V^{1}_{h}]^{n}\times V^{\mathbf{c}}_{h} \times V^{\mathbf{d}}_{h}$ such that $\boldsymbol{U}_{h}= \mathcal{I}^{1}_{h}(\overline{\boldsymbol{U}})$, on $\Gamma_1$ and }}
\begin{subequations}\label{CSFEM23}
\begin{alignat}{3}
\llangle\boldsymbol{P}_{h},\mathbf{grad}\,\boldsymbol{\Upsilon}_{h}\rrangle&=\llangle \boldsymbol{B},\boldsymbol{\Upsilon}_{h} \rrangle + \int_{\Gamma_{2}}\overline{\boldsymbol{T}}\boldsymbol{\cdot}\boldsymbol{\Upsilon}_{h} dA, &\quad &\forall \boldsymbol{\Upsilon}_{h}\in [V^{1}_{h,1}]^{n}, \label{CSFEM23_1}\\
\llangle\mathbf{grad}\,\boldsymbol{U}_{h},\boldsymbol{\lambda}_{h}\rrangle - \llangle \boldsymbol{K}_{h},\boldsymbol{\lambda}_{h}\rrangle &= 0,& &\forall \boldsymbol{\lambda}_{h}\in V^{\mathbf{c}}_{h}, \label{CSFEM23_2}\\
\llangle \mathbb{P}(\boldsymbol{K}_{h}),\boldsymbol{\pi}_{h}\rrangle - \llangle \boldsymbol{P}_{h}, \boldsymbol{\pi}_{h} \rrangle  &= 0, & &\forall \boldsymbol{\pi}_{h}\in V^{\mathbf{d}}_{h}. \label{CSFEM23_3}
\end{alignat}
\end{subequations}
\end{minipage}
\bigskip

\begin{rem} As mentioned earlier in Remark \ref{NotWellDefined}, generally speaking, the response function $\mathbb{P}$ is not well-defined as a mapping $H^{\mathbf{c}}(\mathcal{B})\rightarrow H^{\mathbf{d}}(\mathcal{B})$. By considering simple piecewise polynomial deformation gradients, it is easy to see that $\mathbb{P}$ is not necessarily well-defined as a mapping $V^{\mathbf{c}}_{h}\rightarrow V^{\mathbf{d}}_{h}$ as well; For example, see \citep[Section 4]{An2018}. Thus, in general, we have $\mathbb{P}(\boldsymbol{K}_{h})\notin V^{\mathbf{d}}_{h}$ in \eqref{CSFEM23_3}. This equation simply defines the approximate stress $\boldsymbol{P}_{h}$ as the unique $L^{2}$-orthogonal projection of $\mathbb{P}(\boldsymbol{K}_{h})$ on $V^{\mathbf{d}}_{h}$. In the mixed methods introduced in \citep[Section 3.2]{AnFSYa2017}, $\boldsymbol{P}_{h}\in V^{\mathbf{d}}_{h}$ is the $L^{2}$-orthogonal projection of $\mathbb{P}(\boldsymbol{K}_{h})$ on $V^{\mathbf{c}}_{h}$, which means that $\boldsymbol{P}_{h}\in V^{\mathbf{c}}_{h}\cap V^{\mathbf{d}}_{h}$, and therefore, unlike members of $V^{\mathbf{d}}_{h}$ which can be discontinuous along internal faces of $\mathcal{B}_{h}$, $\boldsymbol{P}_{h}$ is forced to be continuous on $\mathcal{B}_{h}$. The implicit assumption $\boldsymbol{P}_{h}\in V^{\mathbf{c}}_{h}$ is unphysical and severely restricts the solution space of $\boldsymbol{P}_{h}$. As a consequence, it was observed that the extension of the finite element method of \citep{AnFSYa2017} to the $3$D case is very challenging.       
\end{rem}

\begin{rem} In \eqref{CSFEM23_2}, the approximate displacement gradient $\boldsymbol{K}_{h}\in V^{\mathbf{c}}_{h}$ is defined as the unique $L^{2}$-orthogonal projection of $\mathbf{grad}\,\boldsymbol{U}_{h}$ on $V^{\mathbf{c}}_{h}$, with $\boldsymbol{K}_{h}\neq\mathbf{grad}\,\boldsymbol{U}_{h}$, in general. The relation between the complex \eqref{gcdTenHilb3D} and the de Rham complex allows to discretize \eqref{gcdTenHilb3D} by using FEEC. In particular, the discrete de Rham complexes introduced in \citep[Section 5.1]{Arnold2010} implies that \eqref{gcdTenHilb3D} can be discretized as 
\begin{equation*}\label{gcdDis3D}
\scalebox{1}{\xymatrix@C=3.5ex{0 \ar[r] & [V^{1}_{h}]^{3} \ar[r]^-{\mathbf{grad}} &V^{\mathbf{c}}_{h} \ar[r]^-{\mathbf{curl}} & V^{\mathbf{d}}_{h} \ar[r]^-{\mathbf{div}}  &[V_{h}]^{3} \ar[r]  &0,  } } 
\end{equation*}
where the finite element spaces $(V^{1}_{h},V^{\mathbf{c}}_{h},V^{\mathbf{d}}_{h},V_{h})$ are associated to one of the following choices of finite elements:
\begin{alignat*}{5}
\big(\mathrm{LE}_{i},&\,\mathrm{NED}^{2}_{i-1}&&,\mathrm{BDM}_{i-2}&,\,\mathrm{DE}_{i-3}\big), ~ i\geq3, \\
\big(\mathrm{LE}_{i},&\,\mathrm{NED}^{2}_{i-1}&&,\mathrm{RT}_{i-1}&,\,\mathrm{DE}_{i-2}\big), ~ i\geq2,\\
\big(\mathrm{LE}_{i},&\,\mathrm{NED}^{1}_{i}&&,\mathrm{BDM}_{i-1}&,\,\mathrm{DE}_{i-2}\big),~ i\geq2,\\
\big(\mathrm{LE}_{i},&\,\mathrm{NED}^{1}_{i}&&,\mathrm{RT}_{i}&,\,\mathrm{DE}_{i-1}\big), ~ i\geq1,
\end{alignat*}
where $\mathrm{LE}_i$ is the Lagrange element of degree $i$, $\mathrm{NED}^{j}_{i}$ is the $i$-th degree N\'{e}d\'{e}lec element of the $j$-th kind \citep{Ned1986}, $\mathrm{BDM}_{i}$ is the Brezzi-Douglas-Marini element of degree $i$ \citep{BDM1985}, $\mathrm{RT}_{i}$ is the Raviart-Thomas element of degree $i$ \citep{RT1977}, and $\mathrm{DE}_i$ is the discontinuous element of degree $i$, see Figure \ref{FEO12}. Thus, if the finite element spaces $(V^{1}_{h},V^{\mathbf{c}}_{h})$ are induced by $(\mathrm{LE}_{i},\mathrm{NED}^{1}_{i})$ or $(\mathrm{LE}_{i+1},\mathrm{NED}^{2}_{i})$, $i\geq1$, the mapping $\mathbf{grad}: [V^{1}_{h}]^{n}\rightarrow V^{\mathbf{c}}_{h}$ will be well-defined. For these choices of finite element spaces, we have $\boldsymbol{K}_{h}=\mathbf{grad}\,\boldsymbol{U}_{h}$ as $\mathbf{grad}\,\boldsymbol{K}_{h}\in V^{\mathbf{c}}_{h}$ for $\boldsymbol{K}_{h}\in[V^{1}_{h}]^{n}$, i.e. the projection of $\mathbf{grad}\,\boldsymbol{U}_{h}\in V^{\mathbf{c}}_{h}$ on $V^{\mathbf{c}}_{h}$ is equal to itself.  
\end{rem}

\section{Stability Analysis}\label{Sec_Conv}
We employ the general theory introduced in \citep{PoRa1994,CaRa1997} for the Galerkin approximation of nonlinear problems to study the convergence of solutions of \eqref{CSFEM23} to regular solutions of the problem \eqref{3DElasMixF}. This theory is summarized in the Appendix. In particular, we write a sufficient inf-sup condition and two other weaker inf-sup conditions. The former condition is a necessary and sufficient condition for the uniqueness of solutions of the linearization of \eqref{CSFEM23}. We mention a computational framework for studying these inf-sup conditions as well and rigorously show that certain choices of finite elements violate these inf-sup conditions. The following analysis is not valid for singular solutions, which may be studied based on the general approximation framework of \citep{BrRaRa1981}.

\subsection{A Sufficient Stability Condition}
For simplicity and without loss of generality, we assume that $\overline{\boldsymbol{U}}=0$ in \eqref{3DElasMixF}. To apply the theory of \citep{PoRa1994,CaRa1997}, we write the problem \eqref{3DElasMixF} in the abstract form \eqref{AbsProb} as follows: Let $Z=[H^{1}_{1}(\mathcal{B})]^{n}\times H^{\mathbf{c}}(\mathcal{B}) \times H^{\mathbf{d}}(\mathcal{B})$, and let $u,y,z\in Z$, where $u=(\boldsymbol{U},\boldsymbol{K},\boldsymbol{P})$, $y=(\boldsymbol{\Upsilon},\boldsymbol{\lambda},\boldsymbol{\pi})$, $z=(\boldsymbol{V},\boldsymbol{M},\boldsymbol{Q})$. 
Then, \eqref{3DElasMixF} can be stated as: Find $u\in Z$ such that
\begin{equation}\label{Abs_MixProblem}
\begin{alignedat}{3}
\langle H(u), y \rangle &= \llangle\boldsymbol{P},\mathbf{grad}\,\boldsymbol{\Upsilon}\rrangle + \llangle\mathbf{grad}\,\boldsymbol{U},\boldsymbol{\lambda}\rrangle - \llangle \boldsymbol{K},\boldsymbol{\lambda}\rrangle\\
 &+ \llangle \mathbb{P}(\boldsymbol{K}),\boldsymbol{\pi}\rrangle - \llangle \boldsymbol{P}, \boldsymbol{\pi} \rrangle  - \llangle \boldsymbol{B},\boldsymbol{\Upsilon} \rrangle - \int_{\Gamma_{2}}\overline{\boldsymbol{T}}\boldsymbol{\cdot}\boldsymbol{\Upsilon} dA = 0, \quad \forall y \in Z. 
\end{alignedat}
\end{equation}
To write an inf-sup condition for the stability of approximations of the above problem, we consider the bilinear form 
\begin{equation}\label{BiLinCSFEM}
b(z, y) = \llangle\boldsymbol{Q},\mathbf{grad}\,\boldsymbol{\Upsilon}\rrangle + \llangle\mathbf{grad}\,\boldsymbol{V},\boldsymbol{\lambda}\rrangle - \llangle \boldsymbol{M},\boldsymbol{\lambda}\rrangle + \llangle \mathsf{A}(\boldsymbol{K})\boldsymbol{:}\boldsymbol{M},\boldsymbol{\pi}\rrangle - \llangle \boldsymbol{Q}, \boldsymbol{\pi} \rrangle, \quad \forall z,y \in Z. 
\end{equation}  
where $\mathsf{A}(\boldsymbol{K})$ is the elasticity tensor in terms of the displacement gradient and $(\mathsf{A}(\boldsymbol{K})\boldsymbol{:}\boldsymbol{M})^{IJ}:=A^{IJRS}M^{RS}$. This bilinear form is the derivative of the mapping $H$ in \eqref{Abs_MixProblem}, see \eqref{BiLinF}.

Suppose $Z_{h}:=[V^{1}_{h,1}]^{n}\times V^{\mathbf{c}}_{h} \times V^{\mathbf{d}}_{h}$, where the finite element spaces $V^{1}_{h,1}$, $V^{\mathbf{c}}_{h}$, and $V^{\mathbf{d}}_{h}$ were introduced in Section \ref{Sec_CSFEM}. Since $Z_{h}$ is a $Z$-conformal finite element space with the approximability property, the conditions (i) and (ii) of the abstract theory of the Appendix are satisfied and therefore, close to a regular solution $u=(\boldsymbol{U},\boldsymbol{K},\boldsymbol{P})$ of \eqref{3DElasMixF}, the discrete problem \eqref{CSFEM23} has a unique solution $u_{h}=(\boldsymbol{U}_{h},\boldsymbol{K}_{h},\boldsymbol{P}_{h})$ that converges to $u$ as $h\rightarrow 0$ if the inf-sup condition \eqref{infsupCond} holds, that is, if there exists a mesh-independent number $\beta>0$ such that 
\begin{equation}\label{CSFEM-infsup}
\underset{y_{h}\in Z_{h}}{\inf}\, \underset{z_{h}\in Z_{h}}{\sup} \frac{b(z_{h},y_{h})}{\|z_{h}\|_{Z} \|y_{h}\|_{Z} } \geq \beta>0,
\end{equation}
where $y_{h}=(\boldsymbol{\Upsilon}_{h},\boldsymbol{\lambda}_{h},\boldsymbol{\pi}_{h})\in Z_{h}$, $z_{h}=(\boldsymbol{V}_{h},\boldsymbol{M}_{h},\boldsymbol{Q}_{h})\in Z_{h}$, the bilinear form $b(z_{h},y_{h})$ is given in \eqref{BiLinCSFEM}, and $\|z_{h}\|^{2}_{Z}=\|\boldsymbol{V}_{h}\|^{2}_{1}+\|\boldsymbol{M}_{h}\|^{2}_{\mathbf{c}}+\|\boldsymbol{Q}_{h}\|^{2}_{\mathbf{d}}$, with 
\begin{align*}
\|\boldsymbol{V}_{h}\|^{2}_{1} &= \llangle \boldsymbol{V}_{h}, \boldsymbol{V}_{h} \rrangle + \llangle \mathbf{grad}\,\boldsymbol{V}_{h}, \mathbf{grad}\,\boldsymbol{V}_{h} \rrangle,\\
\|\boldsymbol{M}_{h}\|^{2}_{\mathbf{c}} &= \llangle \boldsymbol{M}_{h}, \boldsymbol{M}_{h} \rrangle + \llangle \mathbf{curl}\,\boldsymbol{M}_{h}, \mathbf{curl}\,\boldsymbol{M}_{h} \rrangle,\\
\|\boldsymbol{Q}_{h}\|^{2}_{\mathbf{d}} &= \llangle \boldsymbol{Q}_{h}, \boldsymbol{Q}_{h} \rrangle + \llangle \mathbf{div}\,\boldsymbol{Q}_{h}, \mathbf{div}\,\boldsymbol{Q}_{h} \rrangle.
\end{align*}
Notice that \eqref{CSFEM-infsup} depends on the material properties. If the abstract inf-sup condition \eqref{infsupCond} of the Appendix holds, then the discrete linear system \eqref{DisLinProb} has a unique solution for any given data. Using the bilinear form \eqref{BiLinCSFEM}, this linear system reads: Given $\boldsymbol{f}^{1}$, $\boldsymbol{f}^{\mathbf{c}}$, and $\boldsymbol{f}^{\mathbf{d}}$ of $L^{2}$-class, find $(\boldsymbol{Y}_{\!\!h},\boldsymbol{M}_{h},\boldsymbol{Q}_{h})\in Z_{h}$ such that 
\begin{subequations}\label{LinCSFEM23}
\begin{alignat}{3}
\llangle\boldsymbol{Q}_{h},\mathbf{grad}\,\boldsymbol{\Upsilon}_{h}\rrangle&= \llangle \boldsymbol{f}^{1},\boldsymbol{\Upsilon}_{h}\rrangle, &\quad &\forall \boldsymbol{\Upsilon}_{h}\in [V^{1}_{h,1}]^{n}, \label{LinCSFEM23_I}\\
\llangle\mathbf{grad}\,\boldsymbol{Y}_{\!\!h},\boldsymbol{\lambda}_{h}\rrangle - \llangle \boldsymbol{M}_{h},\boldsymbol{\lambda}_{h}\rrangle &= \llangle \boldsymbol{f}^{\mathbf{c}},\boldsymbol{\lambda}_{h}\rrangle, & &\forall \boldsymbol{\lambda}_{h}\in V^{\mathbf{c}}_{h}, \label{LinCSFEM23_II}\\
\llangle \mathsf{A}(\boldsymbol{K})\boldsymbol{:}\boldsymbol{M}_{h},\boldsymbol{\pi}_{h}\rrangle - \llangle \boldsymbol{Q}_{h}, \boldsymbol{\pi}_{h} \rrangle  &= \llangle \boldsymbol{f}^{\mathbf{d}},\boldsymbol{\pi}_{h}\rrangle, & &\forall \boldsymbol{\pi}_{h}\in V^{\mathbf{d}}_{h}. \label{LinCSFEM23_III}
\end{alignat}
\end{subequations}
Thus, if the material-dependent inf-sup condition \eqref{CSFEM-infsup} holds, the linear system \eqref{LinCSFEM23} will have a unique solution for any input data $\boldsymbol{f}^{1}$, $\boldsymbol{f}^{\mathbf{c}}$, and $\boldsymbol{f}^{\mathbf{d}}$.   
   
Following the computational framework discussed in the Appendix, to computationally investigate the inf-sup condition \eqref{CSFEM-infsup}, we write its matrix form. Given a mesh $\mathcal{B}_{h}$ of the body $\mathcal{B}$, let $\{\boldsymbol{\Psi}_{i}\}_{i=1}^{n_{1}}$, $\{\boldsymbol{\Lambda}_{i}\}_{i=1}^{n_{\mathbf{c}}}$, and $\{\boldsymbol{\Phi}_{i}\}_{i=1}^{n_{\mathbf{d}}}$ be respectively the global shape functions of $[V^{1}_{h,1}]^{n}$, $V^{\mathbf{c}}_{h}$, and $V^{\mathbf{d}}_{h}$ and let $n_{t}=n_{1}+n_{\mathbf{c}}+n_{\mathbf{d}}$ denote the total number of degrees of freedom. Using the relations
\begin{equation*}
\boldsymbol{Y}_{\!\!h}=\sum_{j=1}^{n_{1}}y_{j}\boldsymbol{\Psi}_{j},\quad \boldsymbol{M}_{h}=\sum_{j=1}^{n_{\mathbf{c}}}m_{j}\boldsymbol{\Lambda}_{j}, \text{ and } \boldsymbol{Q}_{h}=\sum_{j=1}^{n_{\mathbf{d}}}q_{j}\boldsymbol{\Phi}_{j},  
\end{equation*}
one can write \eqref{LinCSFEM23} in the matrix form 
\begin{equation}\label{LinSysCSFEM}
\mathbb{S}_{n_{t}\times n_{t}}\cdot\mathbf{z}_{n_{t}\times 1}=\mathbf{f}_{n_{t}\times 1},
\end{equation}
with
\begin{equation}\label{CoefDef}
\mathbb{S}=\left[ \arraycolsep=1.1pt\def\arraystretch{1.2} \begin{array}{c;{2pt/2pt}c;{2pt/2pt}c} \mathbf{0} & \mathbf{0} & \mathbb{S}^{1\mathbf{d}}_{n_{1}\times n_{\mathbf{d}}} \\ \hdashline[2pt/2pt] 
                            \mathbb{S}^{\mathbf{c}1}_{n_{\mathbf{c}}\times n_{1}} & \mathbb{S}^{\mathbf{c}\mathbf{c}}_{n_{\mathbf{c}}\times n_{\mathbf{c}}} &\mathbf{0} \\ \hdashline[2pt/2pt]
														\mathbf{0} & \mathbb{S}^{\mathbf{d}\mathbf{c}}_{n_{\mathbf{d}}\times n_{\mathbf{c}}} & \mathbb{S}^{\mathbf{d}\mathbf{d}}_{n_{\mathbf{d}}\times n_{\mathbf{d}}}\end{array}   \right], \quad
\mathbf{z}=\left[ \arraycolsep=1.1pt\def\arraystretch{1.2} \begin{array}{c} \mathbf{y}_{n_{1}\times1} \\ \mathbf{m}_{n_{\mathbf{c}}\times1} \\  \mathbf{q}_{n_\mathbf{d}\times1}	 \end{array} \right], \text{ and }
\mathbf{f}=\left[ \arraycolsep=1.1pt\def\arraystretch{1.2} \begin{array}{c} \mathbf{f}^{1}_{n_{1}\times1} \\  \mathbf{f}^{\mathbf{c}}_{n_{\mathbf{c}}\times1} \\  \mathbf{f}^{\mathbf{d}}_{n_\mathbf{d}\times1}	 \end{array} \right],					
\end{equation}
where the components of the above matrices and vectors are given by
\begin{alignat*}{5}
& \mathbb{S}^{1\mathbf{d}}_{ij}=\llangle \boldsymbol{\Phi}_{j},\mathbf{grad}\,\boldsymbol{\Psi}_{i}\rrangle,  &~~&\mathbb{S}^{\mathbf{c}1}_{ij}=\llangle \mathbf{grad}\,\boldsymbol{\Psi}_{j},\boldsymbol{\Lambda}_{i}\rrangle,  &~~&\mathbb{S}^{\mathbf{c}\mathbf{c}}_{ij}= -\llangle \boldsymbol{\Lambda}_{i},\boldsymbol{\Lambda}_{j}\rrangle,\\
& \mathbb{S}^{\mathbf{d}\mathbf{c}}_{ij}= \llangle \mathsf{A}(\boldsymbol{K})\boldsymbol{:}\boldsymbol{\Lambda}_{j},\boldsymbol{\Phi}_{i}\rrangle,  && \mathbb{S}^{\mathbf{d}\mathbf{d}}_{ij}= -\llangle \boldsymbol{\Phi}_{i},\boldsymbol{\Phi}_{j}\rrangle, && \\
& \mathbf{y}_{i}= y_{i}, && \mathbf{m}_{i}=m_{i}, &&\mathbf{q}_{i}=q_{i},  \\
&\mathbf{f}^{1}_{i}=\llangle \mathbf{f}^{1}, \boldsymbol{\Psi}_{i} \rrangle, && \mathbf{f}^{\mathbf{c}}_{i}=\llangle \mathbf{f}^{\mathbf{c}}, \boldsymbol{\Lambda}_{i} \rrangle, &&\mathbf{f}^{\mathbf{d}}_{i}= \llangle \mathbf{f}^{\mathbf{d}}, \boldsymbol{\Phi}_{i} \rrangle. 
\end{alignat*}
By replacing the matrix $\mathbb{B}$ of the inf-sup condition \eqref{infsupCondMat} with $\mathbb{S}$, one obtains the matrix form of \eqref{CSFEM-infsup} as 
\begin{equation}\label{infsupCondMat_CSFEM}
\underset{\mathbf{w}\in\mathbb{R}^{n_{t}}}{\inf}\, \underset{\mathbf{z}\in\mathbb{R}^{n_{t}}}{\sup} \frac{\mathbf{w}^{T}\mathbb{S}\, \mathbf{z}}{\|\mathbf{w}\|_{Z} \|\mathbf{z}\|_{Z}} \geq \beta>0,
\end{equation}
with $\|\mathbf{z}\|^{2}_{Z}=\mathbf{z}^{T}\mathbb{D}\,\mathbf{z}$, where the symmetric, positive definite matrix $\mathbb{D}$ is given by
\begin{equation*}
\mathbb{D}_{n_{t}\times n_{t}}=\left[ \arraycolsep=1.1pt\def\arraystretch{1.2} \begin{array}{c;{2pt/2pt}c;{2pt/2pt}c} \mathbb{D}^{1}_{n_{1}\times n_{1}} &\mathbf{0} & \mathbf{0}  \\ \hdashline[2pt/2pt] 
                            \mathbf{0} & \mathbb{D}^{\mathbf{c}}_{n_{\mathbf{c}}\times n_{\mathbf{c}}} &\mathbf{0} \\ \hdashline[2pt/2pt]
														\mathbf{0} & \mathbf{0} & \mathbb{D}^{\mathbf{d}}_{n_{\mathbf{d}}\times n_{\mathbf{d}}}\end{array}   \right],
\end{equation*}
and the components of the symmetric, positive definite matrices $\mathbb{D}^{1}$, $\mathbb{D}^{\mathbf{c}}$, and $\mathbb{D}^{\mathbf{d}}$ are   
\begin{equation}\label{MetMat}
\begin{aligned}
\mathbb{D}^{1}_{ij} &= \llangle \boldsymbol{\Psi}_{i}, \boldsymbol{\Psi}_{j} \rrangle + \llangle \mathbf{grad}\,\boldsymbol{\Psi}_{i}, \mathbf{grad}\,\boldsymbol{\Psi}_{j} \rrangle,\\
\mathbb{D}^{\mathbf{c}}_{ij} &= \llangle \boldsymbol{\Lambda}_{i}, \boldsymbol{\Lambda}_{j} \rrangle + \llangle \mathbf{curl}\,\boldsymbol{\Lambda}_{i}, \mathbf{curl}\,\boldsymbol{\Lambda}_{j} \rrangle,\\
\mathbb{D}^{\mathbf{d}}_{ij} &= \llangle \boldsymbol{\Phi}_{i}, \boldsymbol{\Phi}_{j} \rrangle + \llangle \mathbf{div}\,\boldsymbol{\Phi}_{i}, \mathbf{div}\,\boldsymbol{\Phi}_{j} \rrangle.
\end{aligned}
\end{equation}
The discussion of the Appendix then implies that the inf-sup condition \eqref{CSFEM-infsup} holds if and only if the smallest singular value of $\mathbb{M}^{Z}\mathbb{S}\,\mathbb{M}^{Z}$ is positive and bounded from below by a positive number $\beta$ as $h\rightarrow 0$, where $\mathbb{M}^{Z}$ is the unique symmetric, positive definite matrix that satisfies $(\mathbb{M}^{Z})^{2}=\mathbb{D}$.

\subsection{Weaker Stability Conditions}
If the inf-sup condition \eqref{CSFEM-infsup} holds, then \eqref{LinCSFEM23} will have a unique solution for any input data. In particular, \eqref{LinCSFEM23_I} must have a solution for any $\boldsymbol{f}^{1}$, or equivalently, the left-hand side of \eqref{LinCSFEM23_I} must define an onto mapping. This condition is equivalent to the following material-independent inf-sup condition: There exists $\alpha_{h}>0$ such that
\begin{equation}\label{CSFEM-infsup_Sur}
\underset{\boldsymbol{\Upsilon}_{h}\in [V^{1}_{h,1}]^{n}}{\inf}\, \underset{\boldsymbol{Q}_{h}\in V^{\mathbf{d}}_{h}}{\sup} \frac{\llangle\boldsymbol{Q}_{h},\mathbf{grad}\,\boldsymbol{\Upsilon}_{h}\rrangle}{\|\boldsymbol{Q}_{h}\|_{\mathbf{d}}\, \|\boldsymbol{\Upsilon}_{h}\|_{1} } \geq \alpha_{h}.
\end{equation}
On the other hand, \eqref{LinCSFEM23_I} and \eqref{LinCSFEM23_III} must have a solution for any $\boldsymbol{f}^{1}$ and $\boldsymbol{f}^{\mathbf{d}}$, which means that the left-hand side of \eqref{LinCSFEM23_I} and \eqref{LinCSFEM23_III} should define an onto mapping. This latter condition can be stated by another inf-sup condition: There should be $\gamma_{h}>0$ such that 
\begin{equation}\label{CSFEM-infsup_Sur2}
\underset{(\boldsymbol{\Upsilon}_{h},\boldsymbol{\pi}_{h})\in Z_{1\mathbf{d}}}{\inf}\, \underset{(\boldsymbol{M}_{h},\boldsymbol{Q}_{h})\in Z_{\mathbf{cd}}}{\sup} \frac{\llangle\boldsymbol{Q}_{h},\mathbf{grad}\,\boldsymbol{\Upsilon}_{h}\rrangle + \llangle \mathsf{A}(\boldsymbol{K})\boldsymbol{:}\boldsymbol{M}_{h},\boldsymbol{\pi}_{h}\rrangle - \llangle \boldsymbol{Q}_{h}, \boldsymbol{\pi}_{h} \rrangle  }{\|(\boldsymbol{M}_{h},\boldsymbol{Q}_{h})\|_{\mathbf{cd}} \, \|(\boldsymbol{\Upsilon}_{h},\boldsymbol{\pi}_{h})\|_{1\mathbf{d}} } \geq \gamma_{h},
\end{equation}
with $\|(\boldsymbol{M}_{h},\boldsymbol{Q}_{h})\|^{2}_{\mathbf{cd}} = \|\boldsymbol{M}_{h}\|^{2}_{\mathbf{c}} + \|\boldsymbol{Q}_{h}\|^{2}_{\mathbf{d}}$, and $\|(\boldsymbol{\Upsilon}_{h},\boldsymbol{\pi}_{h})\|^{2}_{\mathbf{1d}} = \|\boldsymbol{\Upsilon}_{h}\|^{2}_{1} + \|\boldsymbol{\pi}_{h}\|^{2}_{\mathbf{d}}$. The inf-sup conditions \eqref{CSFEM-infsup_Sur} and \eqref{CSFEM-infsup_Sur2} are weaker than \eqref{CSFEM-infsup} in the sense that they are only necessary for the validity of \eqref{CSFEM-infsup}.

The material-independent inf-sup condition \eqref{CSFEM-infsup_Sur} admits the matrix form
\begin{equation}\label{CSFEM-infsup_SurMat}
\underset{\mathbf{y}\in\mathbb{R}^{n_{1}}}{\inf}\, \underset{\mathbf{q}\in\mathbb{R}^{n_{\mathbf{d}}}}{\sup} \frac{\mathbf{y}^{T}\mathbb{S}^{1\mathbf{d}} \mathbf{q}}{\|\mathbf{y}\|_{1} \|\mathbf{q}\|_{\mathbf{d}}} \geq \alpha_{h}>0,
\end{equation}
where $\mathbb{S}^{1\mathbf{d}}_{n_{1}\times n_{\mathbf{d}}}$ is defined in \eqref{CoefDef}. Let $\mathbb{M}^{1}$ and $\mathbb{M}^{\mathbf{d}}$ be the unique symmetric and positive definite matrices such that $(\mathbb{M}^{1})^{2}=\mathbb{D}^{1}$, and $(\mathbb{M}^{\mathbf{d}})^{2}=\mathbb{D}^{\mathbf{d}}$, where $\mathbb{D}^1$ and $\mathbb{D}^{\mathbf{d}}$ are given in \eqref{MetMat}. Then, \eqref{CSFEM-infsup_Sur} holds if and only if the smallest singular value of $\mathbb{M}^{1}\mathbb{S}^{1\mathbf{d}}\mathbb{M}^{\mathbf{d}}$ is positive. Similarly, the matrix form of \eqref{CSFEM-infsup_Sur2} reads
\begin{equation}\label{CSFEM-infsup_SurMat2}
\underset{\mathbf{u}\in\mathbb{R}^{n_{1}+n_{\mathbf{d}}}}{\inf}\, \underset{\mathbf{x}\in\mathbb{R}^{n_{\mathbf{c}}+n_{\mathbf{d}}}}{\sup} \frac{\mathbf{u}^{T}\mathbb{G} \mathbf{x}}{\|\mathbf{u}\|_{1\mathbf{d}} \|\mathbf{x}\|_{\mathbf{cd}}} \geq \gamma_{h}>0,
\end{equation}
with 
\begin{equation*}
\mathbb{G}=\left[ \arraycolsep=1.1pt\def\arraystretch{1.2} \begin{array}{c;{2pt/2pt}c} \mathbf{0} & \mathbb{S}^{1\mathbf{d}}_{n_{1}\times n_{\mathbf{d}}} \\ \hdashline[2pt/2pt]
														\mathbb{S}^{\mathbf{d}\mathbf{c}}_{n_{\mathbf{d}}\times n_{\mathbf{c}}} & \mathbb{S}^{\mathbf{d}\mathbf{d}}_{n_{\mathbf{d}}\times n_{\mathbf{d}}}\end{array}   \right], 			
\end{equation*}
where $\mathbb{S}^{1\mathbf{d}}$, $\mathbb{S}^{\mathbf{dc}}$, and $\mathbb{S}^{\mathbf{dd}}$ are defined in \eqref{CoefDef}. Suppose $\mathbb{M}^{1\mathbf{d}}$ and $\mathbb{M}^{\mathbf{cd}}$ are positive definite matrices that satisfy
\begin{equation*}
(\mathbb{M}^{1\mathbf{d}})^{2}=\left[ \arraycolsep=1.1pt\def\arraystretch{1.2} \begin{array}{c;{2pt/2pt}c} \mathbb{D}^{1}_{n_{1}\times n_{1}} & \mathbf{0} \\ \hdashline[2pt/2pt]
\mathbf{0} & \mathbb{D}^{\mathbf{d}}_{n_{\mathbf{d}}\times n_{\mathbf{d}}} \end{array}   \right], \text{ and } (\mathbb{M}^{\mathbf{cd}})^{2}=\left[ \arraycolsep=1.1pt\def\arraystretch{1.2} \begin{array}{c;{2pt/2pt}c} \mathbb{D}^{\mathbf{c}}_{n_{\mathbf{c}}\times n_{\mathbf{c}}} & \mathbf{0} \\ \hdashline[2pt/2pt]
\mathbf{0} & \mathbb{D}^{\mathbf{d}}_{n_{\mathbf{d}}\times n_{\mathbf{d}}} \end{array}   \right], 			
\end{equation*}
where $\mathbb{D}^{1}$, $\mathbb{D}^{\mathbf{c}}$, and $\mathbb{D}^{\mathbf{d}}$ were introduced in \eqref{MetMat}. The condition \eqref{CSFEM-infsup_Sur2} holds if and only if the smallest singular value of $\mathbb{M}^{1\mathbf{d}}\mathbb{G}\,\mathbb{M}^{\mathbf{cd}}$ is positive.

The inf-sup condition \eqref{CSFEM-infsup_SurMat} is equivalent to the surjectivity of the linear mapping $\mathbb{S}^{1\mathbf{d}}:\mathbb{R}^{n_{\mathbf{d}}}\rightarrow\mathbb{R}^{n_{1}}$, that is, the matrix $\mathbb{S}^{1\mathbf{d}}$ being full ranked. This result can be directly deduced from the structure of the matrix $\mathbb{S}$ in \eqref{LinSysCSFEM} as well. As a consequence of the rank-nullity theorem, it is easy to see that $\mathbb{S}^{1\mathbf{d}}$ is not full rank if $n_{\mathbf{d}}<n_{1}$. The upshot can be stated as follows.    

\begin{thm}\label{SimplInfSup} Suppose $n_{1} = \dim([V^{1}_{h,1}]^{n})$, and $n_{\mathbf{d}} = \dim V^{\mathbf{d}}_{h}$. The inf-sup conditions \eqref{CSFEM-infsup_Sur} and \eqref{CSFEM-infsup} do not hold if $n_{\mathbf{d}}<n_{1}$.
\end{thm}

The condition \eqref{CSFEM-infsup_SurMat2} is equivalent to the surjectivity of $\mathbb{G}: \mathbb{R}^{n_{\mathbf{c}}+ n_{\mathbf{d}}} \rightarrow \mathbb{R}^{n_{\mathbf{1}}+ n_{\mathbf{d}}}$, that is, $\mathbb{G}$ being full rank. Due to the rank-nullity theorem, this result does not hold if $n_{\mathbf{c}}+ n_{\mathbf{d}} < n_{1}+ n_{\mathbf{d}}$. Thus, one concludes that: 

\begin{thm}\label{SimplInfSup2} Suppose $n_{1} = \dim([V^{1}_{h,1}]^{n})$, and $n_{\mathbf{c}} = \dim V^{\mathbf{c}}_{h}$. The inf-sup conditions \eqref{CSFEM-infsup_Sur2} and \eqref{CSFEM-infsup} do not hold if $n_{\mathbf{c}}<n_{1}$.
\end{thm}

Notice that if the inf-sup condition \eqref{CSFEM-infsup_Sur} fails, then the discrete nonlinear problem \eqref{CSFEM23} is not stable as it may not have any solution for some body forces and boundary tractions. Assume that $\mathcal{B}$ is a polyhedral domain with a triangular (2D) or a tetrahedral (3D) mesh $\mathcal{B}_{h}$, which is geometrically conformal. Let $N_{v}$, $N_{ed}$, and $N_{f}$ be respectively the number of vertices, edges, and faces of $\mathcal{B}_{h}$ (in $2$D, we have $N_{f}=N_{ed}$). For the $n$-dimensional elements $\text{LE}_{2}$, $\text{NED}^{1}_{1}$, and $\text{RT}_{1}$, $n=2,3$, of Figure \ref{FEO12}, it is straightforward to show that $n_{1}= n(N_{v}+N_{ed})$, $n_{\mathbf{c}} = n\,N_{ed}$, and $n_{\mathbf{d}} = n\,N_{f}$. These relations imply the following corollary of Theorems \ref{SimplInfSup} and \ref{SimplInfSup2}.

\begin{cor}\label{Cor_infsup} Let $\text{FE}_{\mathbf{c}}$ and $\text{FE}_{\mathbf{d}}$ respectively be arbitrary $H^{\mathbf{c}}$- and $H^{\mathbf{d}}$-conformal finite elements. We have:
\begin{enumerate}
\item In $2$D, the finite element choice $(\text{LE}_{2},\text{FE}_{\mathbf{c}},\text{RT}_{1})$ for mixed finite element methods \eqref{CSFEM23} does not satisfy the inf-sup conditions \eqref{CSFEM-infsup_Sur} and \eqref{CSFEM-infsup}.
\item In $2$D and $3$D, the finite element choice $(\text{LE}_{2},\text{NED}^{1}_{1},\text{FE}_{\mathbf{d}})$ for mixed finite element methods \eqref{CSFEM23} does not satisfy the inf-sup conditions \eqref{CSFEM-infsup_Sur2} and \eqref{CSFEM-infsup}.
\end{enumerate} 
\end{cor}

\begin{rem} In \citep{Re2015}, a three-field formulation for linearized elasticity in terms of displacement, strain, and stress was introduced, which is similar to the system \eqref{LinCSFEM23} but by using discontinuous $L^{2}$-elements instead of $H^{\mathbf{c}}$- and $H^{\mathbf{d}}$-conformal elements. For that linear system, it was shown that the ellipticity of the elasticity tensor and the analogue of the inf-sup condition \eqref{CSFEM-infsup_Sur} in terms of $L^{2}$ finite element spaces are sufficient for the well-posedness \citep[Theorem 5.2]{Re2015}. The inf-sup condition \eqref{CSFEM-infsup} is a stronger condition in the sense that it is both necessary and sufficient for the well-posedness of \eqref{LinCSFEM23}.    
\end{rem}

\begin{rem}The condition \eqref{CSFEM-infsup_Sur} is similar to the Babu\v{s}ka-Brezzi condition for the Stokes problem. For choices of finite elements that \eqref{CSFEM-infsup_Sur} fails, one can use strategies similar to those for the Babu\v{s}ka-Brezzi condition to enrich $V^{\mathbf{d}}_{h}$, e.g. employing bubble functions or using a finer mesh for $V^{\mathbf{d}}_{h}$, see \citep[Chapter 4]{ErnGuermond2004}.        
\end{rem}

\section{Numerical Results}\label{Sec_Examples}
To study the performance of the mixed finite element method \eqref{CSFEM23}, we employ the finite elements of Figure \ref{FEO12} and solve several $2$D and $3$D numerical examples in this section. Numerical simulations are performed by using FEniCS \citep{loMaWe2012}, which is an open-source platform with the high-level Python and C++ interfaces. For our simulations, we consider compressible Neo-Hookean materials with the stored energy function
\begin{equation*}
W(\boldsymbol{F})=\frac{\mu}{2}(I_{1}-3) - \frac{\mu}{2}\ln I_{3} + \frac{\lambda}{2}(\ln I_{3})^{2}, \quad \mu,\lambda>0,
\end{equation*}
where $\boldsymbol{F}$ is the deformation gradient, $I_{1}=\mathrm{tr}\,\boldsymbol{C}$, and $I_{3}=\det \boldsymbol{C}$, with $\boldsymbol{C}=\boldsymbol{F}^{T}\boldsymbol{F}$. The constitutive equation in terms of $\boldsymbol{F}$ then reads
\begin{equation*}\label{NeoHook} 
\mathbb{P}(\boldsymbol{F})= \mu\boldsymbol{F} -\mu \boldsymbol{F}^{-T} + 2\lambda (\ln I_{3}) \boldsymbol{F}^{-T}.
\end{equation*}
Substituting $\boldsymbol{F}=\boldsymbol{I}+\boldsymbol{K}$ in the above equation yields the constitutive equation $\mathbb{P}(\boldsymbol{K})$ in terms of the displacement gradient $\boldsymbol{K}$, where $\boldsymbol{I}$ is the identity matrix. Moreover, the tensor $\mathsf{A}(\boldsymbol{K})\boldsymbol{:}\boldsymbol{M}$ in the bilinear form \eqref{BiLinCSFEM} becomes
\begin{equation*}
\mathsf{A}(\boldsymbol{K})\boldsymbol{:}\boldsymbol{M}= \mu \boldsymbol{M} + (\mu-2\lambda\ln I_{3}) \boldsymbol{F}^{-T}\boldsymbol{M}^{T}\boldsymbol{F}^{-T} + 4\lambda(\mathrm{tr}\, \boldsymbol{F}^{-1}\boldsymbol{M}) \boldsymbol{F}^{-T}.
\end{equation*}

\paragraph{Convention}To concisely refer to a choice of the elements of Figure \ref{FEO12} for the mixed finite element methods \eqref{CSFEM23}, we use the following convention: $\mathrm{L}i$, $\mathrm{N}ji$, $\mathrm{R}i$, and $\mathrm{B}i$ respectively denote $\text{LE}_i$, $\text{NED}^{j}_{i}$, $\text{RT}_{i}$, and $\mathrm{BDM}_{i}$. For example, $\mathrm{L1N12B2}$ denotes the choice $(\text{LE}_{1},\text{NED}^{1}_{2},\text{BDM}_{2})$.

\begin{figure}[thb]
\begin{center}
\includegraphics[scale=.52,angle=0]{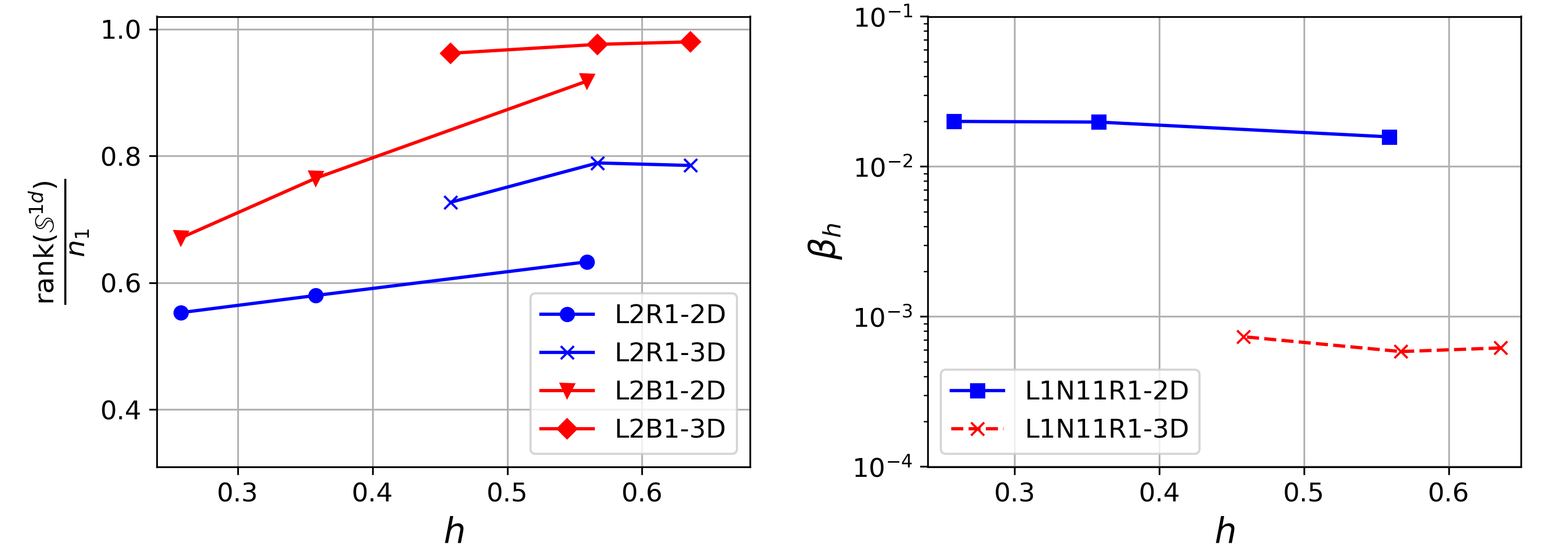}
\end{center}
\vspace*{-0.3in}
\caption{\footnotesize Numerical analysis of the inf-sup conditions \eqref{CSFEM-infsup_Sur} (the left panel) and \eqref{CSFEM-infsup} (the right panel) using unstructured meshes of the unit square (2D) and the unit cube (3D) similar to the second row of Figures \ref{PlateStUnst} and \ref{CubeStUnstM}. Left Panel: The ratio $\mathrm{rank}(\mathbb{S}^{1\mathbf{d}})/n_{1}$ versus the maximum diameter of elements $h$ for the finite element choices $\mathrm{L2R1}$ and $\mathrm{L2B1}$ for $(\boldsymbol{U},\boldsymbol{P})$. As this ratio is smaller than $1$, $\mathbb{S}^{1\mathbf{d}}$ is rank deficient and the condition \eqref{CSFEM-infsup_Sur} does not hold. Right Panel: Values of the lower bound $\beta_h$ of the inf-sup condition \eqref{CSFEM-infsup} versus $h$ associated to $2$D and $3$D meshes and the choice of elements $\mathrm{L1N11R1}$ for $(\boldsymbol{U},\boldsymbol{K},\boldsymbol{P})$. The results suggest that $\beta_h$ does not decrease as $h$ decreases, and hence there is a positive lower bound $\beta$ that satisfies \eqref{CSFEM-infsup}.} 
\label{InfSupR}
\end{figure}

\begin{figure}[tbh]
\begin{center}
\includegraphics[scale=.4,angle=0]{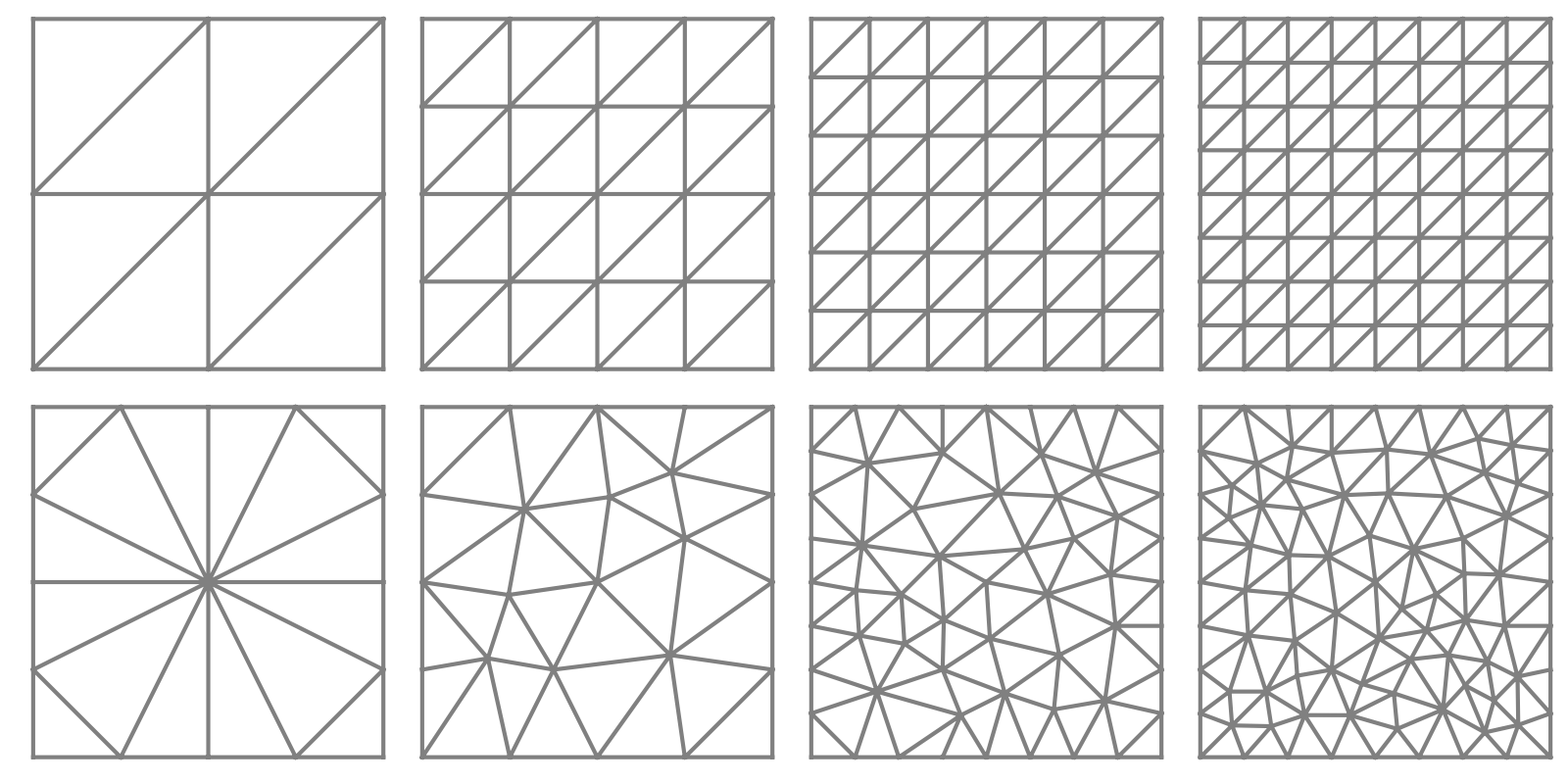}
\end{center}
\vspace*{-0.3in}
\caption{\footnotesize Meshes of a unit square where the number of elements  $N_{e}$ and the maximum diameter of elements $h$ are given by  $(N_{e},h)= (8,0.707)$, $(32,0.354)$, $(72,0.236)$, $(128,0.177)$, for the structured meshes of the first row and $(N_{e},h) = (16, 0.559)$, $(32,0.358)$, $(88,0.249)$, $(146,0.181)$, for the unstructured meshes of the second row.} 
\label{PlateStUnst}
\end{figure}

\subsection{Stability Analysis}\label{subsec_StAn}
We begin by numerically investigating the inf-sup conditions introduced earlier by using their matrix forms. Of course, the following numerical results are not mathematical proofs; Rather, they provide strong evidences for obtaining mathematical proofs.         

As discussed earlier, the inf-sup conditions \eqref{CSFEM-infsup_Sur} and \eqref{CSFEM-infsup_Sur2} are respectively equivalent to $\mathbb{S}^{1\mathbf{d}}$ and $\mathbb{G}$ being full rank. Therefore, to study the validity of these conditions, one can study the rank of the associated matrices. Due to Corollary \ref{Cor_infsup}, we already know that the choice $\mathrm{L2R1}$ for $(\boldsymbol{U},\boldsymbol{P})$ in $2$D does not satisfy the inf-sup condition \eqref{CSFEM-infsup_Sur}. Our numerical studies show that this choice does not satisfy \eqref{CSFEM-infsup_Sur} in $3$D as well. Moreover, the choice $\mathrm{L2B1}$ does not satisfy \eqref{CSFEM-infsup_Sur} neither in $2$D nor in $3$D. For example, the left panel of Figure \ref{InfSupR} depicts the ratio $\mathrm{rank}(\mathbb{S}^{1\mathbf{d}})/n_{1}$ versus the maximum diameter of elements $h$ for the choices $\mathrm{L2R1}$ and $\mathrm{L2B1}$ calculated using unstructured meshes of the unit square and the unit cube. The matrix $\mathbb{S}^{1\mathbf{d}}$ is not full rank in these cases since $\frac{\mathrm{rank}(\mathbb{S}^{1\mathbf{d}})}{n_{1}} <1$. It is interesting to note that unlike $\mathrm{L2R1}$, we have $n_{\mathbf{d}}> n_{1}$ for $\mathrm{L2B1}$.

Corollary \ref{Cor_infsup} implies that the choice $\mathrm{L2N11}$ for $(\boldsymbol{U},\boldsymbol{K})$ does not satisfy the inf-sup condition \eqref{CSFEM-infsup_Sur2}. Notice that unlike \eqref{CSFEM-infsup_Sur}, the inf-sup condition  \eqref{CSFEM-infsup_Sur2} is material-dependent. Our numerical experiments suggest that for Neo-Hookean materials with regular deformations, all other choices of the elements of Figure \ref{FEO12} for $(\boldsymbol{U},\boldsymbol{K})$ satisfy \eqref{CSFEM-infsup_Sur2}.     

The inf-sup condition \eqref{CSFEM-infsup} is sufficient for the convergence of solutions of \eqref{CSFEM23} to regular solutions of \eqref{3DElasMixF}. Our numerical results for Neo-Hookean materials discussed in the remainder of this section suggest that all choices of elements of Figure \ref{FEO12} for $(\boldsymbol{U},\boldsymbol{K},\boldsymbol{P})$ that satisfy the inf-sup conditions \eqref{CSFEM-infsup_Sur} and \eqref{CSFEM-infsup_Sur2} also satisfy the inf-sup condition \eqref{CSFEM-infsup}. As an example, the right panel of Figure \ref{InfSupR} shows the values of the lower bound $\beta_h$ of the matrix form \eqref{infsupCondMat_CSFEM} for the choice of elements $\mathrm{L1N11R1}$ in $2$D and $3$D. Results are calculated using unstructured meshes of the unit square and the unit cube with the material parameters $\mu=\lambda =1$ near the reference configuration, that is, $\boldsymbol{K}=0$. To approximate $\beta_h$ for each mesh, one can employ the smallest singular value of $\mathbb{M}^{Z}\mathbb{S}\,\mathbb{M}^{Z}$ or equivalently, the square root of the smallest eigenvalue of $\mathbb{D}\,\mathbb{S}^{T}\mathbb{D}\,\mathbb{S}$, with $\mathbb{D} = (\mathbb{M}^{Z})^{2}$. The results suggests that the values of $\beta_h$ are bounded from below as $h$ decreases and therefore, there is a lower bound $\beta>0$ that satisfies \eqref{CSFEM-infsup}.  

The validity of the material-dependent inf-sup conditions \eqref{CSFEM-infsup} and \eqref{CSFEM-infsup_Sur2} is dependent to properties of the elasticity tensor $\mathsf{A}(\boldsymbol{K})$. To rigorously study these inf-sup conditions, one needs to impose some additional restrictions on $\mathsf{A}$ which are physically reasonable. The classical inequalities for $\mathsf{A}$ \citep[Section 51]{TrNo65} and suitable assumptions on the stored energy functional $W$ such as polyconvexity \citep{Ci1988} are relevant here.

\begin{table}[!thb]
\tabcolsep=3.0pt
\caption{Convergence rates $r$ and $L^{2}$-errors of the plate example: DoF is the number of total degrees of freedom and $(E_{\boldsymbol{U}}, E_{\boldsymbol{K}}, E_{\boldsymbol{P}})=(\|\boldsymbol{U}_{h}-\boldsymbol{U}_{e}\|, \|\boldsymbol{K}_{h}-\boldsymbol{K}_{e}\|, \|\boldsymbol{P}_{h}-\boldsymbol{P}_{e}\|)$ are the $L^{2}$-errors of the approximate solution ($\boldsymbol{U}_{h}$, $\boldsymbol{F}_{h}$, $\boldsymbol{P}_{h}$) with respect to the exact solution ($\boldsymbol{U}_{e}$, $\boldsymbol{F}_{e}$, $\boldsymbol{P}_{e}$) associated to \eqref{2Dplate_Exact}.}
\resizebox{0.8\textwidth}{!}{\begin{minipage}{\textwidth}
\centering
\begin{tabular}{cl|cc|cc|cc||cl|cc|cc|cc}
\toprule
FEM & DoF &$E_{\boldsymbol{U}}$ &  & $E_{\boldsymbol{K}}$ &  & $E_{\boldsymbol{P}}$ &  & FEM & DoF &$E_{\boldsymbol{U}}$ &  & $E_{\boldsymbol{K}}$ &  & $E_{\boldsymbol{P}}$ &  \\
\midrule
\multirow{4}{*}{\rotatebox[origin=c]{90}{\parbox[c]{1.5cm}{\centering $\mathrm{L1N11R1}$}}}
& 82   & 1.76e-2 & \multirow{4}{*}{\rotatebox[origin=c]{90}{\parbox[c]{1.5cm}{\centering $r=2.0$}}}  & 1.21e-1 & \multirow{4}{*}{\rotatebox[origin=c]{90}{\parbox[c]{1.5cm}{\centering $r=1.0$}}}  & 1.66e-1 & \multirow{4}{*}{\rotatebox[origin=c]{90}{\parbox[c]{1.5cm}{\centering $r=1.0$}}} & \multirow{4}{*}{\rotatebox[origin=c]{90}{\parbox[c]{1.5cm}{\centering $\mathrm{L1N11R2}$}}} & 146   & 1.82e-2 & \multirow{4}{*}{\rotatebox[origin=c]{90}{\parbox[c]{1.5cm}{\centering $r=2.1$}}}  & 1.20e-1 & \multirow{4}{*}{\rotatebox[origin=c]{90}{\parbox[c]{1.5cm}{\centering $r=1.0$}}}  & 1.61e-1 &\multirow{4}{*}{\rotatebox[origin=c]{90}{\parbox[c]{1.5cm}{\centering $r=1.0$}}} \\
& 274   & 4.26e-3 &  & 6.49e-2 &  & 9.22e-2 & & & 514   & 4.51e-3 &  & 6.24e-2 &  & 7.98e-2 &\\
& 578  & 2.04e-3 &  & 4.76e-2 &  & 9.14e-2 & & & 1106   & 1.95e-3 &  & 4.19e-2 &  & 5.23e-2 &\\
& 994  & 1.07e-3 &  & 3.30e-2 &  & 4.66e-2 & & & 1922   & 1.08e-3 &  & 3.15e-2 &  & 3.87e-2 &\\
\bottomrule
\multirow{4}{*}{\rotatebox[origin=c]{90}{\parbox[c]{1.5cm}{\centering $\mathrm{L1N12R1}$}}}
& 146   & 1.75e-2 & \multirow{4}{*}{\rotatebox[origin=c]{90}{\parbox[c]{1.5cm}{\centering $r=2.0$}}}  & 1.21e-1 & \multirow{4}{*}{\rotatebox[origin=c]{90}{\parbox[c]{1.5cm}{\centering $r=1.0$}}}  & 1.66e-1 & \multirow{4}{*}{\rotatebox[origin=c]{90}{\parbox[c]{1.5cm}{\centering $r=1.0$}}} & \multirow{4}{*}{\rotatebox[origin=c]{90}{\parbox[c]{1.5cm}{\centering $\mathrm{L1N12R2}$}}} & 210   & 1.82e-2 & \multirow{4}{*}{\rotatebox[origin=c]{90}{\parbox[c]{1.5cm}{\centering $r=2.1$}}}  & 1.20e-1 & \multirow{4}{*}{\rotatebox[origin=c]{90}{\parbox[c]{1.5cm}{\centering $r=1.0$}}}  & 1.61e-1 &\multirow{4}{*}{\rotatebox[origin=c]{90}{\parbox[c]{1.5cm}{\centering $r=1.0$}}} \\
& 514   & 4.26e-3 &  & 6.49e-2 &  & 9.22e-2 & & & 754   & 4.51e-3 &  & 6.24e-2 &  & 7.98e-2 &\\
& 1106  & 2.04e-3 &  & 4.76e-2 &  & 9.14e-2 & & & 1634   & 1.95e-3 &  & 4.19e-2 &  & 5.23e-2 &\\
& 1922  & 1.07e-3 &  & 3.30e-2 &  & 4.66e-2 & & & 2850   & 1.08e-3 &  & 3.15e-2 &  & 3.87e-2 &\\
\bottomrule
\multirow{4}{*}{\rotatebox[origin=c]{90}{\parbox[c]{1.5cm}{\centering $\mathrm{L1N11B1}$}}}
& 114   & 1.76e-2 & \multirow{4}{*}{\rotatebox[origin=c]{90}{\parbox[c]{1.5cm}{\centering $r=2.1$}}}  & 1.21e-1 & \multirow{4}{*}{\rotatebox[origin=c]{90}{\parbox[c]{1.5cm}{\centering $r=1.0$}}}  & 1.55e-1 & \multirow{4}{*}{\rotatebox[origin=c]{90}{\parbox[c]{1.5cm}{\centering $r=1.1$}}} & \multirow{4}{*}{\rotatebox[origin=c]{90}{\parbox[c]{1.5cm}{\centering $\mathrm{L1N11B2}$}}} & 194   & 1.83e-2 & \multirow{4}{*}{\rotatebox[origin=c]{90}{\parbox[c]{1.5cm}{\centering $r=2.1$}}}  & 1.20e-1 & \multirow{4}{*}{\rotatebox[origin=c]{90}{\parbox[c]{1.5cm}{\centering $r=1.0$}}}  & 1.62e-1 &\multirow{4}{*}{\rotatebox[origin=c]{90}{\parbox[c]{1.5cm}{\centering $r=1.0$}}} \\
& 386   & 4.25e-3 &  & 6.32e-2 &  & 7.63e-2 & & & 690   & 4.55e-3 &  & 6.23e-2 &  & 8.05e-2 &\\
& 818  & 1.84e-3 &  & 4.24e-2 &  & 4.90e-2 & & & 1490   & 1.97e-3 &  & 4.18e-2 &  & 5.28e-2 &\\
& 1410  & 1.02e-3 &  & 3.19e-2 &  & 3.65e-2 & & & 2594   & 1.09e-3 &  & 3.15e-2 &  & 3.91e-2 &\\
\bottomrule
\multirow{4}{*}{\rotatebox[origin=c]{90}{\parbox[c]{1.5cm}{\centering $\mathrm{L1N12B1}$}}}
& 178   & 1.76e-2 & \multirow{4}{*}{\rotatebox[origin=c]{90}{\parbox[c]{1.5cm}{\centering $r=2.1$}}}  & 1.21e-1 & \multirow{4}{*}{\rotatebox[origin=c]{90}{\parbox[c]{1.5cm}{\centering $r=1.0$}}}  & 1.55e-1 & \multirow{4}{*}{\rotatebox[origin=c]{90}{\parbox[c]{1.5cm}{\centering $r=1.1$}}} & \multirow{4}{*}{\rotatebox[origin=c]{90}{\parbox[c]{1.5cm}{\centering $\mathrm{L1N12B2}$}}} & 258   & 1.83e-2 & \multirow{4}{*}{\rotatebox[origin=c]{90}{\parbox[c]{1.5cm}{\centering $r=2.1$}}}  & 1.20e-1 & \multirow{4}{*}{\rotatebox[origin=c]{90}{\parbox[c]{1.5cm}{\centering $r=1.0$}}}  & 1.62e-1 &\multirow{4}{*}{\rotatebox[origin=c]{90}{\parbox[c]{1.5cm}{\centering $r=1.0$}}} \\
& 626   & 4.25e-3 &  & 6.32e-2 &  & 7.63e-2 & & & 930   & 4.55e-3 &  & 6.23e-2 &  & 8.05e-2 &\\
& 1346  & 1.84e-3 &  & 4.24e-2 &  & 4.90e-2 & & & 2018   & 1.97e-3 &  & 4.18e-2 &  & 5.28e-2 &\\
& 2338  & 1.02e-3 &  & 3.19e-2 &  & 3.65e-2 & & & 3522   & 1.09e-3 &  & 3.15e-2 &  & 3.91e-2 &\\
\bottomrule
\multirow{4}{*}{\rotatebox[origin=c]{90}{\parbox[c]{1.5cm}{\centering $\mathrm{L1N21R1}$}}}
& 114   & 1.76e-2 & \multirow{4}{*}{\rotatebox[origin=c]{90}{\parbox[c]{1.5cm}{\centering $r=2.0$}}}  & 1.21e-1 & \multirow{4}{*}{\rotatebox[origin=c]{90}{\parbox[c]{1.5cm}{\centering $r=1.0$}}}  & 1.66e-1 & \multirow{4}{*}{\rotatebox[origin=c]{90}{\parbox[c]{1.5cm}{\centering $r=1.0$}}} & \multirow{4}{*}{\rotatebox[origin=c]{90}{\parbox[c]{1.5cm}{\centering $\mathrm{L1N21R2}$}}} & 178   & 1.82e-2 & \multirow{4}{*}{\rotatebox[origin=c]{90}{\parbox[c]{1.5cm}{\centering $r=2.1$}}}  & 1.20e-1 & \multirow{4}{*}{\rotatebox[origin=c]{90}{\parbox[c]{1.5cm}{\centering $r=1.0$}}}  & 1.61e-1 &\multirow{4}{*}{\rotatebox[origin=c]{90}{\parbox[c]{1.5cm}{\centering $r=1.0$}}} \\
& 386   & 4.26e-3 &  & 6.49e-2 &  & 9.22e-2 & & & 626   & 4.51e-3 &  & 6.24e-2 &  & 7.98e-2 &\\
& 818  & 2.04e-3 &  & 4.76e-2 &  & 9.14e-2 & & & 1346   & 1.95e-3 &  & 4.19e-2 &  & 5.23e-2 &\\
& 1410  & 1.07e-3 &  & 3.30e-2 &  & 4.66e-2 & & & 2338   & 1.08e-3 &  & 3.15e-2 &  & 3.87e-2 &\\
\bottomrule
\multirow{4}{*}{\rotatebox[origin=c]{90}{\parbox[c]{1.5cm}{\centering $\mathrm{L1N22R1}$}}}
& 194   & 1.76e-2 & \multirow{4}{*}{\rotatebox[origin=c]{90}{\parbox[c]{1.5cm}{\centering $r=2.0$}}}  & 1.21e-1 & \multirow{4}{*}{\rotatebox[origin=c]{90}{\parbox[c]{1.5cm}{\centering $r=1.0$}}}  & 1.66e-1 & \multirow{4}{*}{\rotatebox[origin=c]{90}{\parbox[c]{1.5cm}{\centering $r=1.0$}}} & \multirow{4}{*}{\rotatebox[origin=c]{90}{\parbox[c]{1.5cm}{\centering $\mathrm{L1N22R2}$}}} & 258   & 1.82e-2 & \multirow{4}{*}{\rotatebox[origin=c]{90}{\parbox[c]{1.5cm}{\centering $r=2.1$}}}  & 1.20e-1 & \multirow{4}{*}{\rotatebox[origin=c]{90}{\parbox[c]{1.5cm}{\centering $r=1.0$}}}  & 1.61e-1 &\multirow{4}{*}{\rotatebox[origin=c]{90}{\parbox[c]{1.5cm}{\centering $r=1.0$}}} \\
& 690   & 4.26e-3 &  & 6.49e-2 &  & 9.22e-2 & & & 930   & 4.51e-3 &  & 6.24e-2 &  & 7.98e-2 &\\
& 1490  & 2.04e-3 &  & 4.76e-2 &  & 9.14e-2 & & & 2018   & 1.95e-3 &  & 4.19e-2 &  & 5.23e-2 &\\
& 2594  & 1.07e-3 &  & 3.30e-2 &  & 4.66e-2 & & & 3522   & 1.08e-3 &  & 3.15e-2 &  & 3.87e-2 &\\
\bottomrule
\multirow{4}{*}{\rotatebox[origin=c]{90}{\parbox[c]{1.5cm}{\centering $\mathrm{L1N21B1}$}}}
& 146   & 1.76e-2 & \multirow{4}{*}{\rotatebox[origin=c]{90}{\parbox[c]{1.5cm}{\centering $r=2.1$}}}  & 1.21e-1 & \multirow{4}{*}{\rotatebox[origin=c]{90}{\parbox[c]{1.5cm}{\centering $r=1.0$}}}  & 1.55e-1 & \multirow{4}{*}{\rotatebox[origin=c]{90}{\parbox[c]{1.5cm}{\centering $r=1.1$}}} & \multirow{4}{*}{\rotatebox[origin=c]{90}{\parbox[c]{1.5cm}{\centering $\mathrm{L1N21B2}$}}} & 226   & 1.83e-2 & \multirow{4}{*}{\rotatebox[origin=c]{90}{\parbox[c]{1.5cm}{\centering $r=2.1$}}}  & 1.20e-1 & \multirow{4}{*}{\rotatebox[origin=c]{90}{\parbox[c]{1.5cm}{\centering $r=1.0$}}}  & 1.62e-1 &\multirow{4}{*}{\rotatebox[origin=c]{90}{\parbox[c]{1.5cm}{\centering $r=1.0$}}} \\
& 498   & 4.25e-3 &  & 6.32e-2 &  & 7.63e-2 & & & 802   & 4.55e-3 &  & 6.23e-2 &  & 8.05e-2 &\\
& 1058  & 1.84e-3 &  & 4.24e-2 &  & 4.90e-2 & & & 1730   & 1.97e-3 &  & 4.18e-2 &  & 5.28e-2 &\\
& 1826  & 1.02e-3 &  & 3.19e-2 &  & 3.65e-2 & & & 3010   & 1.09e-3 &  & 3.15e-2 &  & 3.91e-2 &\\
\bottomrule
\multirow{4}{*}{\rotatebox[origin=c]{90}{\parbox[c]{1.5cm}{\centering $\mathrm{L1N22B1}$}}}
& 226   & 1.76e-2 & \multirow{4}{*}{\rotatebox[origin=c]{90}{\parbox[c]{1.5cm}{\centering $r=2.1$}}}  & 1.21e-1 & \multirow{4}{*}{\rotatebox[origin=c]{90}{\parbox[c]{1.5cm}{\centering $r=1.0$}}}  & 1.55e-1 & \multirow{4}{*}{\rotatebox[origin=c]{90}{\parbox[c]{1.5cm}{\centering $r=1.1$}}} & \multirow{4}{*}{\rotatebox[origin=c]{90}{\parbox[c]{1.5cm}{\centering $\mathrm{L1N22B2}$}}} & 306   & 1.83e-2 & \multirow{4}{*}{\rotatebox[origin=c]{90}{\parbox[c]{1.5cm}{\centering $r=2.1$}}}  & 1.20e-1 & \multirow{4}{*}{\rotatebox[origin=c]{90}{\parbox[c]{1.5cm}{\centering $r=1.0$}}}  & 1.62e-1 &\multirow{4}{*}{\rotatebox[origin=c]{90}{\parbox[c]{1.5cm}{\centering $r=1.0$}}} \\
& 802   & 4.25e-3 &  & 6.32e-2 &  & 7.63e-2 & & & 1106   & 4.55e-3 &  & 6.23e-2 &  & 8.05e-2 &\\
& 1730  & 1.84e-3 &  & 4.24e-2 &  & 4.90e-2 & & & 2402   & 1.97e-3 &  & 4.18e-2 &  & 5.28e-2 &\\
& 3010  & 1.02e-3 &  & 3.19e-2 &  & 3.65e-2 & & & 4194   & 1.09e-3 &  & 3.15e-2 &  & 3.91e-2 &\\
\bottomrule
\multirow{4}{*}{\rotatebox[origin=c]{90}{\parbox[c]{1.5cm}{\centering $\mathrm{L2N21R2}$}}}
& 210   & 1.30e-3 & \multirow{4}{*}{\rotatebox[origin=c]{90}{\parbox[c]{1.5cm}{\centering $r=2.9$}}}  & 1.71e-2 & \multirow{4}{*}{\rotatebox[origin=c]{90}{\parbox[c]{1.5cm}{\centering $r=1.9$}}}  & 2.40e-2 & \multirow{4}{*}{\rotatebox[origin=c]{90}{\parbox[c]{1.5cm}{\centering $r=1.9$}}} & \multirow{4}{*}{\rotatebox[origin=c]{90}{\parbox[c]{1.5cm}{\centering $\mathrm{L2N21B2}$}}} & 258   & 1.30e-3 & \multirow{4}{*}{\rotatebox[origin=c]{90}{\parbox[c]{1.5cm}{\centering $r=3.0$}}}  & 1.71e-2 & \multirow{4}{*}{\rotatebox[origin=c]{90}{\parbox[c]{1.5cm}{\centering $r=1.9$}}}  & 2.41e-2 &\multirow{4}{*}{\rotatebox[origin=c]{90}{\parbox[c]{1.5cm}{\centering $r=2.0$}}} \\
& 738   & 1.67e-4 &  & 4.44e-3 &  & 6.18e-3 & & & 914   & 1.66e-4 &  & 4.39e-3 &  & 6.18e-3 &\\
& 1586  & 5.05e-5 &  & 2.04e-3 &  & 2.79e-3 & & & 1970   & 5.00e-5 &  & 2.00e-3 &  & 2.78e-3 &\\
& 2754  & 2.19e-5 &  & 1.19e-3 &  & 1.60e-3 & & & 3426   & 2.14e-5 &  & 1.16e-3 &  & 1.59e-3 &\\
\bottomrule
\multirow{4}{*}{\rotatebox[origin=c]{90}{\parbox[c]{1.5cm}{\centering $\mathrm{L2N12R2}$}}}
& 242   & 1.30e-3 & \multirow{4}{*}{\rotatebox[origin=c]{90}{\parbox[c]{1.5cm}{\centering $r=2.9$}}}  & 1.71e-2 & \multirow{4}{*}{\rotatebox[origin=c]{90}{\parbox[c]{1.5cm}{\centering $r=1.9$}}}  & 2.40e-2 & \multirow{4}{*}{\rotatebox[origin=c]{90}{\parbox[c]{1.5cm}{\centering $r=1.9$}}} & \multirow{4}{*}{\rotatebox[origin=c]{90}{\parbox[c]{1.5cm}{\centering $\mathrm{L2N12B2}$}}} & 290   & 1.30e-3 & \multirow{4}{*}{\rotatebox[origin=c]{90}{\parbox[c]{1.5cm}{\centering $r=3.0$}}}  & 1.71e-2 & \multirow{4}{*}{\rotatebox[origin=c]{90}{\parbox[c]{1.5cm}{\centering $r=1.9$}}}  & 2.41e-2 &\multirow{4}{*}{\rotatebox[origin=c]{90}{\parbox[c]{1.5cm}{\centering $r=2.0$}}} \\
& 866   & 1.67e-4 &  & 4.43e-3 &  & 6.18e-3 & & & 1042   & 1.66e-4 &  & 4.39e-3 &  & 6.18e-3 &\\
& 1874  & 5.05e-5 &  & 2.04e-3 &  & 2.79e-3 & & & 2258   & 5.00e-5 &  & 2.00e-3 &  & 2.78e-3 &\\
& 3266  & 2.19e-5 &  & 1.19e-3 &  & 1.60e-3 & & & 3938   & 2.15e-5 &  & 1.16e-3 &  & 1.59e-3 &\\
\bottomrule
\multirow{4}{*}{\rotatebox[origin=c]{90}{\parbox[c]{1.5cm}{\centering $\mathrm{L2N22R2}$}}}
& 290   & 1.30e-3 & \multirow{4}{*}{\rotatebox[origin=c]{90}{\parbox[c]{1.5cm}{\centering $r=2.9$}}}  & 1.71e-2 & \multirow{4}{*}{\rotatebox[origin=c]{90}{\parbox[c]{1.5cm}{\centering $r=1.9$}}}  & 2.40e-2 & \multirow{4}{*}{\rotatebox[origin=c]{90}{\parbox[c]{1.5cm}{\centering $r=1.9$}}} & \multirow{4}{*}{\rotatebox[origin=c]{90}{\parbox[c]{1.5cm}{\centering $\mathrm{L2N22B2}$}}} & 338   & 1.30e-3 & \multirow{4}{*}{\rotatebox[origin=c]{90}{\parbox[c]{1.5cm}{\centering $r=3.0$}}}  & 1.71e-2 & \multirow{4}{*}{\rotatebox[origin=c]{90}{\parbox[c]{1.5cm}{\centering $r=1.9$}}}  & 2.41e-2 &\multirow{4}{*}{\rotatebox[origin=c]{90}{\parbox[c]{1.5cm}{\centering $r=2.0$}}} \\
& 1042   & 1.67e-4 &  & 4.43e-3 &  & 6.18e-3 & & & 1218   & 1.66e-4 &  & 4.39e-3 &  & 6.18e-3 &\\
& 2258  & 5.05e-5 &  & 2.04e-3 &  & 2.79e-3 & & & 2642   & 5.00e-5 &  & 2.00e-3 &  & 2.78e-3 &\\
& 3938  & 2.19e-5 &  & 1.19e-3 &  & 1.60e-3 & & & 4610   & 2.15e-5 &  & 1.16e-3 &  & 1.59e-3 &\\
\bottomrule
\end{tabular} 
\label{Pl1ConvError}
\end{minipage} }
\end{table}

\subsection{Deformation of a $2$D Plate}\label{Ex_2DPlate}
To studying the convergence rate of solutions, we consider a unit-square plate with the material parameters $\mu=\lambda=1$ and solve \eqref{CSFEM23} by employing the body force and the boundary conditions that induce the displacement field 
\begin{equation}\label{2Dplate_Exact}
\boldsymbol{U}_{e} = \left[\begin{array}{c} \frac{1}{2}Y^{3} + \frac{1}{2}\sin(\frac{\pi}{2}Y) \\ 0 \end{array} \right],
\end{equation}
where $(X,Y)$ denotes the Cartesian coordinates in $\mathbb{R}^{2}$. We use Newton's method to solve the resulting nonlinear systems. The linear system solved in each Newton iteration is similar to the linear system \eqref{LinCSFEM23} where $\boldsymbol{K}$ is replaced with the solution of the previous iteration. Therefore, the coefficient matrix of each Newton iteration is similar to the matrix $\mathbb{S}$ of the inf-sup condition \eqref{infsupCondMat_CSFEM} and consequently, Newton iterations become singular if any of the inf-sup conditions introduced earlier (with the solution of the previous iteration instead of $\boldsymbol{K}$) is not satisfied.    

Table \ref{Pl1ConvError} shows $L^{2}$-errors and the associated convergence rates of the mixed method \eqref{CSFEM23} which are calculated by using the structured meshes in the first row of Figure \ref{PlateStUnst} with different combinations of the $2$D elements of degrees 1 and 2 of Figure \ref{FEO12}. The convergence rate $r$ means the error is $O(h^{r})$ as $h\rightarrow 0$, where $h$ is the maximum diameter of elements of a mesh. We observe that $22$ combinations out of $32$ possible combinations of the $2$D elements of Figure \ref{FEO12} are stable. More specifically, the $10$ unstable cases include $\mathrm{L2N11R}i$, $\mathrm{L2N11B}i$, $\mathrm{L2N}ij\mathrm{R1}$, and $\mathrm{L2N}ij\mathrm{B1}$, $i,j=1,2$. Following Corollary \ref{Cor_infsup} and the results of Section \ref{subsec_StAn}, we already know that the cases $\mathrm{L2N}ij\mathrm{R1}$ and $\mathrm{L2N}ij\mathrm{B1}$ are unstable as they do not satisfy the inf-sup condition \eqref{CSFEM-infsup_Sur} and that $\mathrm{L2N11R}i$ and $\mathrm{L2N11B}i$ are unstable as they do not satisfy the inf-sup condition \eqref{CSFEM-infsup_Sur2}. Thus, the inf-sup conditions \eqref{CSFEM-infsup_Sur} and \eqref{CSFEM-infsup_Sur2} are sufficient for studying the stability of this example.

Table \ref{Pl1ConvError} suggests that methods with the element $\mathrm{L1}$ for displacement have very close errors and convergence rates regardless of the degrees of their elements for $\boldsymbol{K}$ and $\boldsymbol{P}$. A similar conclusion also holds for methods with the element $\mathrm{L2}$. This suggests that the degree of the element for displacement has a significant effect on the overall performance of these mixed finite element methods. The optimal convergence rate (that is, the convergence rate of finite element interpolations of sufficiently smooth functions) of $\mathrm{L}i$, $\mathrm{N2i}$ and $\mathrm{B}i$ is $i+1$ while that of $\mathrm{N1i}$ and $\mathrm{R}i$ is $i$ \citep{BoBrFo2013}. Table \ref{Pl1ConvError} shows that the convergence rates for displacement gradient and stress may not be optimal but the convergence rate of displacement is always optimal.

\begin{figure}[tbh]
\begin{center}
\includegraphics[scale=.5,angle=0]{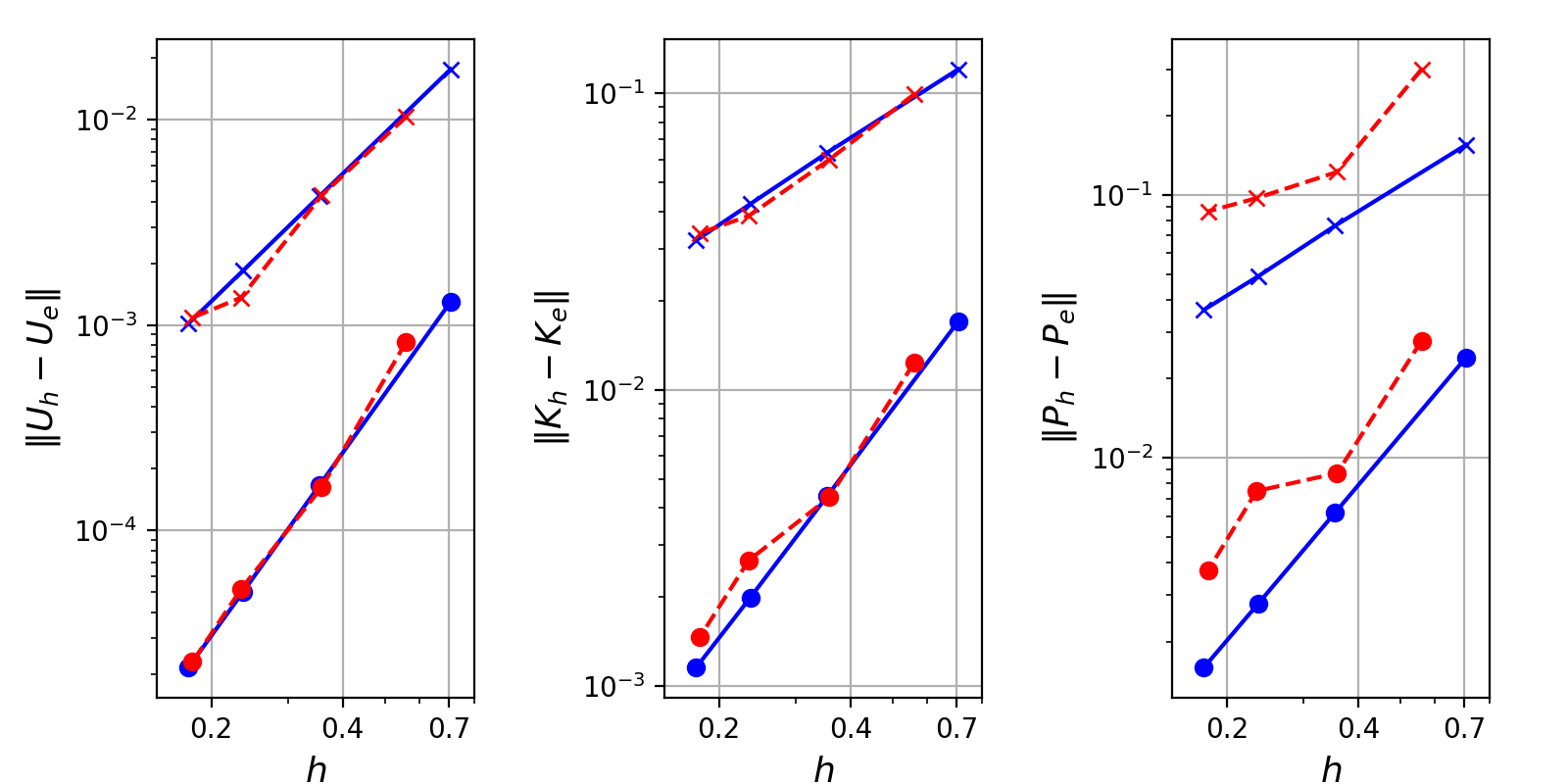}
\end{center}
\vspace*{-0.25in}
\caption{\footnotesize $L^{2}$-errors of displacement $\|\boldsymbol{U}_{h} - \boldsymbol{U}_{e}\|$, displacement gradient $\|\boldsymbol{K}_{h} - \boldsymbol{K}_{e}\|$, and stress $\|\boldsymbol{P}_{h} - \boldsymbol{P}_{e}\|$ associated to the structured meshes (the solid lines) and the unstructured meshes (the dashed lines) of Figure \ref{PlateStUnst}. The data marked by $\times$ and $\bullet$ are respectively calculated by the first-order elements $\mathrm{L1N21B1}$ and the second-order elements $\mathrm{L2N22B2}$.} 
\label{2DConvStUnst}
\end{figure}

\begin{figure}[tbh]
\begin{center}
\includegraphics[scale=.08,angle=0]{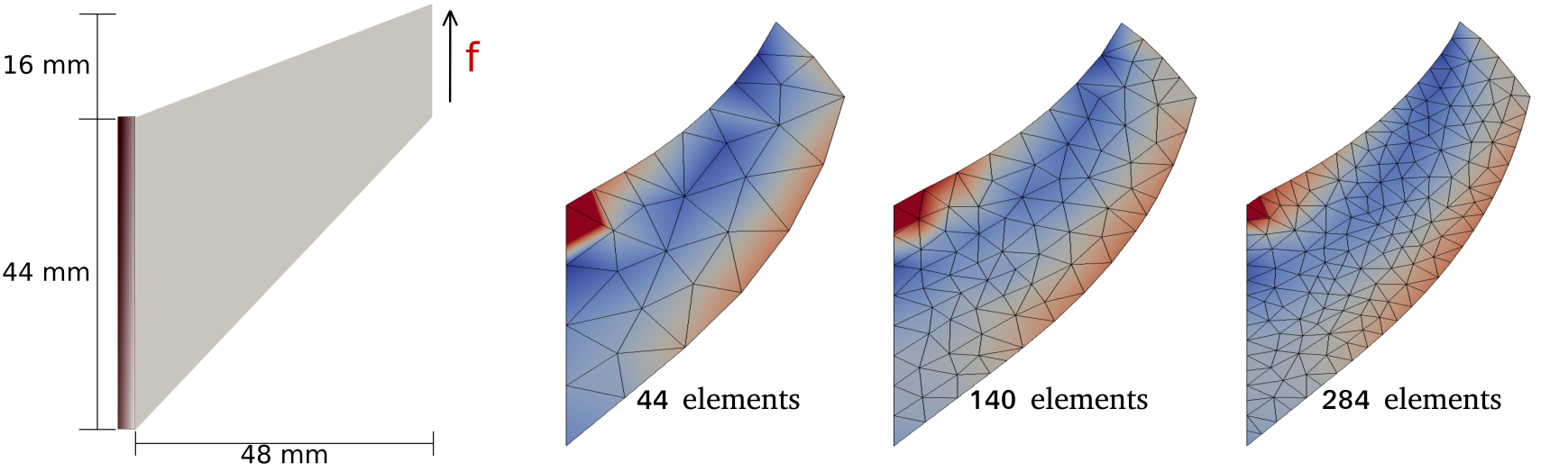}
\end{center}
\vspace*{-0.3in}
\caption{\footnotesize The geometry and deformed configurations of $2$D Cook's membrane. The deformed configurations are computed using the elements $\mathrm{L2N22B2}$ and the shear force $f = 24\,\mathrm{N}/\mathrm{mm}$. Colors in the deformed configuration depict the distribution of the Frobenius norm of stress $\|\boldsymbol{P}\|_{f}$.} 
\label{CookConfig}
\end{figure}

To compare the formulation of this paper with that of \citep{AnFSYa2017}, we notice that the latter mixed formulation is stable only for $7$ out of $32$ possible combinations of the elements of Figure \ref{FEO12}. A comparison between Table \ref{Pl1ConvError} and Table 3 of \citep{AnFSYa2017} suggests that the performance of $\mathrm{L1N11R1}$, $\mathrm{L1N12B1}$, and $\mathrm{L1N22B1}$ is nearly similar in these two formulations while the performance of $\mathrm{L1N12R1}$, $\mathrm{L1N22R1}$, $\mathrm{L1N22R2}$, and $\mathrm{L2N22R2}$ is better using the formulation of this paper. As will be shown in the sequel, the main advantage of the present formulation is that unlike the formulation of \citep{AnFSYa2017} which is only stable in $2$D,  it is stable in both 2D and 3D.

For the brevity, we only consider the choices $\mathrm{L1N21B1}$ and $\mathrm{L2N22B2}$ to solve the other 2D examples of this section. To study the effect of mesh irregularities on the performance of these finite element methods, the $L^{2}$-norm of errors corresponding to structured and unstructured meshes of Figure \ref{PlateStUnst} are shown in Figure \ref{2DConvStUnst}. These results suggest that comparing to the accuracy of approximate displacement and displacement gradient, mesh irregularities may have more impact on the accuracy of approximate stress. Notice that the slope of the curves in Figure \ref{2DConvStUnst} which are associated to the structured meshes are the convergence rates of Table \ref{Pl1ConvError}.

\begin{figure}[tbh]
\begin{center}
\includegraphics[scale=.8,angle=0]{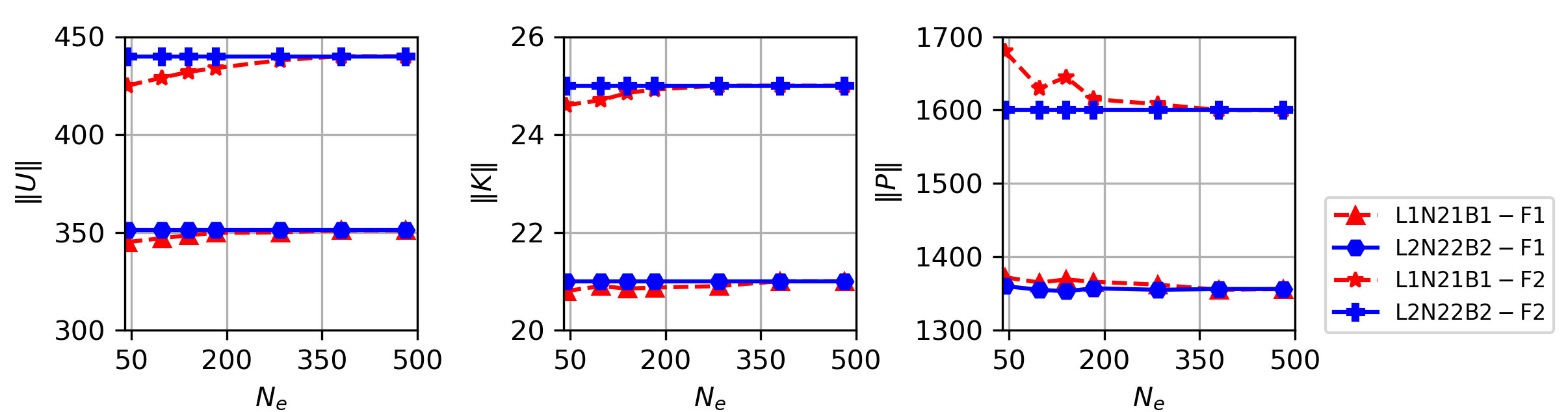}
\end{center}
\vspace*{-0.3in}
\caption{\footnotesize The $L^{2}$-norm of approximate solutions versus the number of elements of meshes $N_{e}$ for $2$D Cook's membrane. The dashed and the solid lines correspond to $\mathrm{L1N21B1}$ and $\mathrm{L2N22B2}$, respectively. Results are computed using unstructured meshes and the shear forces $F1 = 24\,\mathrm{N/mm}$ and $F2 = 32\,\mathrm{N/mm}$.}
\label{L2Norm}
\end{figure}

\begin{figure}[tbh]
\begin{center}
\includegraphics[scale=0.08,angle=0]{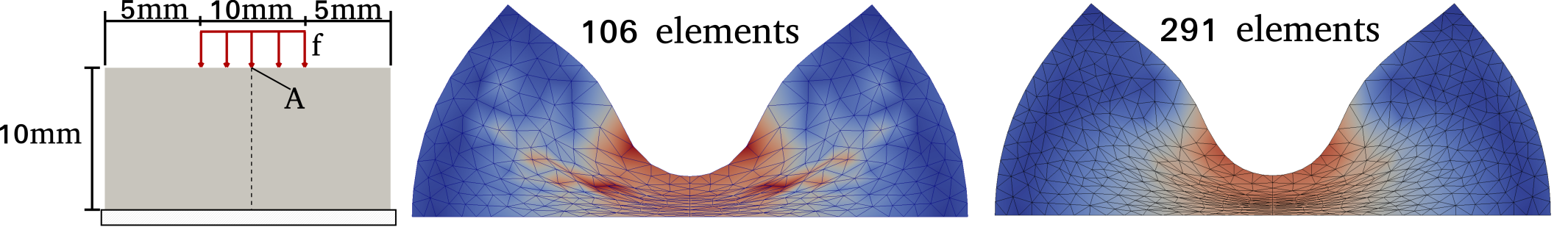}
\end{center}
\vspace*{-0.3in}
\caption{\footnotesize The geometry and deformed configurations of the inhomogeneous compression example. The deformed configurations are computed using the elements $\mathrm{L2N22B2}$ and the force $f = 600\,\mathrm{N}/\mathrm{mm}$. Colors in the deformed configuration depict the distribution of the Frobenius norm of stress $\|\boldsymbol{P}\|_{f}$. }
\label{2D_Inhomo}
\end{figure}

\begin{figure}[tbh]
\begin{center}
\includegraphics[scale=.7,angle=0]{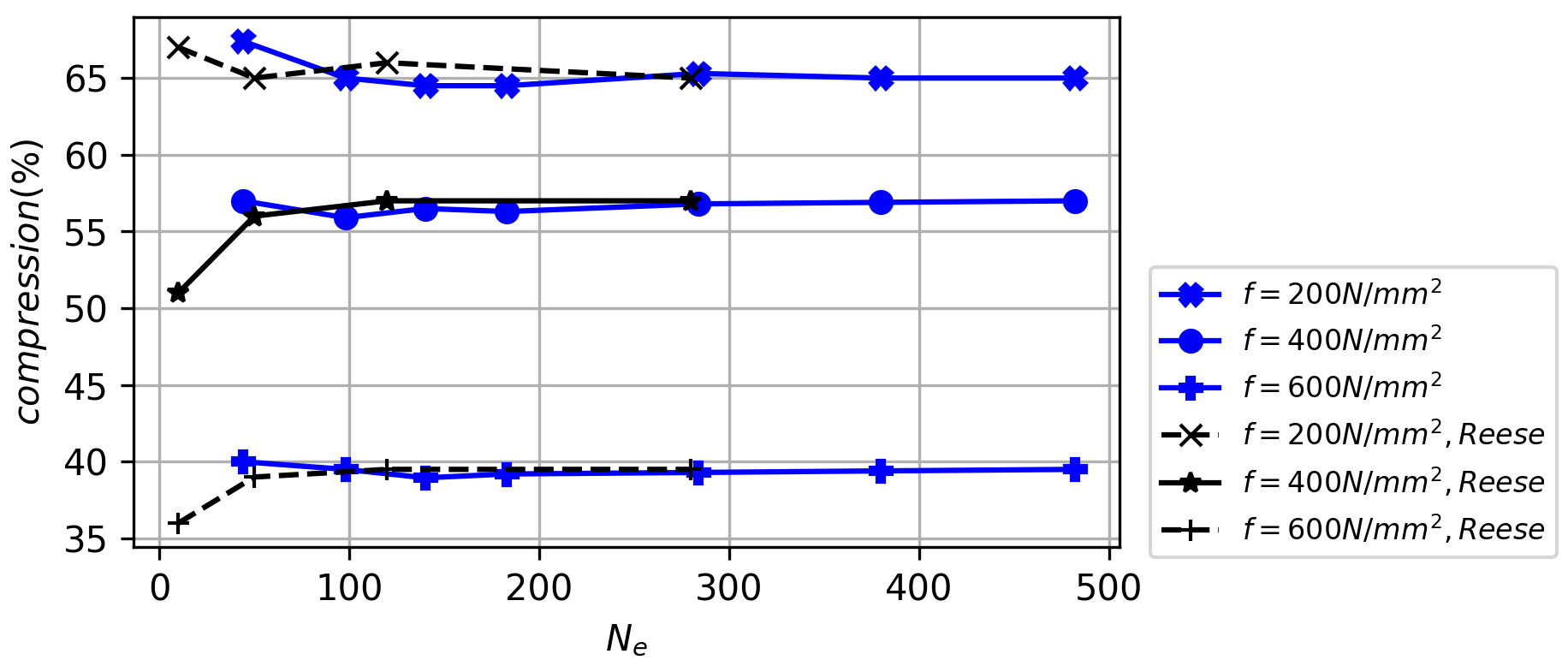}
\end{center}
\vspace*{-0.3in}
\caption{\footnotesize The percentage of the compression of the point $A$ of Figure \ref{2D_Inhomo} versus the number of elements $N_{e}$. Results are computed using the elements $\mathrm{L2N22B2}$. The associated results of \citet{Re2002} are also shown for the comparison.}
\label{Inhomo_comp}
\end{figure}

\subsection{$2$D Cook's Membrane}
Consider the $2$D Cook's membrane problem with the geometry shown in Figure \ref{CookConfig}. This example is usually used to study the performance in bending and in the near-incompressible regime \citep{Re2002}. The material properties are $\mu = 80.194\,\mathrm{N}/\mathrm{mm}^{2}$, and $\lambda = 400889.8\,\mathrm{N}/\mathrm{mm}^{2}$.

Figure \ref{CookConfig} shows deformed configurations calculated using the element $\mathrm{L1N21B1}$ and the load $f = 24\,\mathrm{N}/\mathrm{mm}$. Colors in the deformed configurations depict the distribution of the Frobenius norm of stress $\|\boldsymbol{P}\|_{f} = \sqrt{\mathrm{tr}\,\boldsymbol{P}^{T}\boldsymbol{P}} = \sqrt{\sum_{I,J}|P^{IJ}|^{2}}$. Figure \ref{L2Norm} shows the convergence of the $L^{2}$-norms of approximate solutions. Results are calculated using the elements $\mathrm{L1N21B1}$ and $\mathrm{L2N22B2}$ with two different loads of magnitudes $24$ and $32\,\mathrm{N/mm}$. These results suggest that the mixed formulation \eqref{CSFEM23} can provide accurate approximations of stress in bending and in the near-compressible regime.

\begin{figure}[tbh]
\begin{center}
\includegraphics[scale=.8,angle=0]{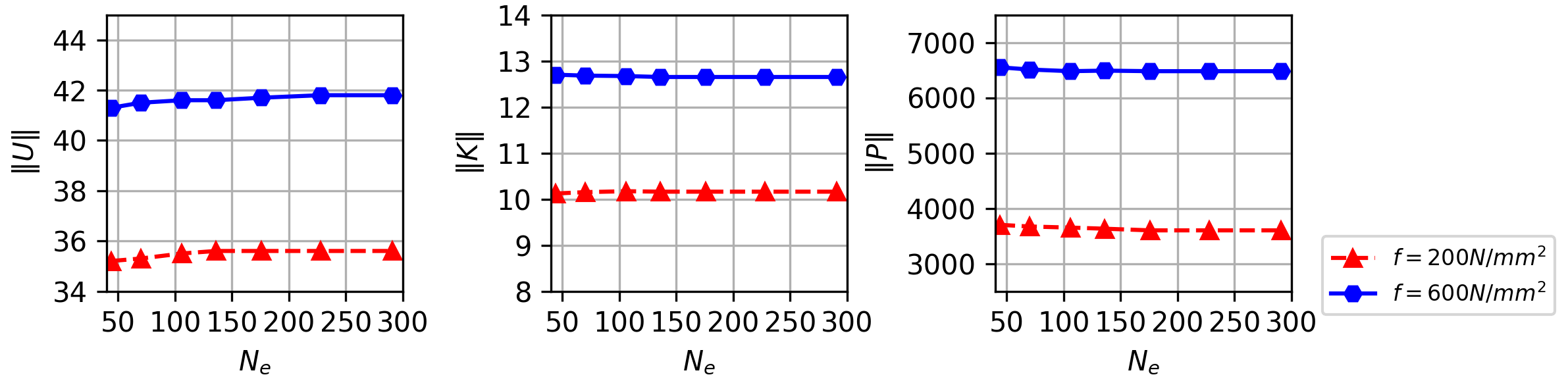}
\end{center}
\vspace*{-0.3in}
\caption{\footnotesize The $L^{2}$-norms of solutions of the $2$D inhomogeneoud compression problem versus the number of elements $N_{e}$. The elements $\mathrm{L2N22B2}$ were used for computing these results.}
\label{InHo_norm}
\end{figure}

\subsection{Inhomogeneous Compression}
Enhanced strain methods are nonconformal three-field methods for small and finite deformations \citep{SiAr1992}. It is well-known that in some cases, these methods may become unstable due to the so-called hourglass instability \citep{Re2002}. One example for this type of instability is the inhomogeneous compression problem shown in Figure \ref{2D_Inhomo}. The horizontal displacement at the top of the domain and the vertical displacement at the bottom are assumed to be zero and the material properties are the same as the previous example.

Deformed configurations of this problem associated to two different meshes which are calculated using the elements $\mathrm{L2N22B2}$  and the load $f = 600\,\mathrm{N}/\mathrm{mm}$ are shown in Figure \ref{2D_Inhomo}. Colors in the deformed configurations show the distribution of the Frobenius norm of stress. Figure \ref{Inhomo_comp} depicts the percentage of compression versus the number of elements for different loads $f$. The compression level is calculated using the vertical displacement of the point $A$ of Figure \ref{2D_Inhomo}, which is located at the midpoint of the top boundary. The results are consistent with those of \citep{Re2002}. Figure \ref{InHo_norm} shows the convergence of the $L^{2}$-norm of solutions by refining meshes. We do not observe any numerical instability in our computations.

\begin{figure}[tbh]
\begin{center}
\includegraphics[scale=.12,angle=0]{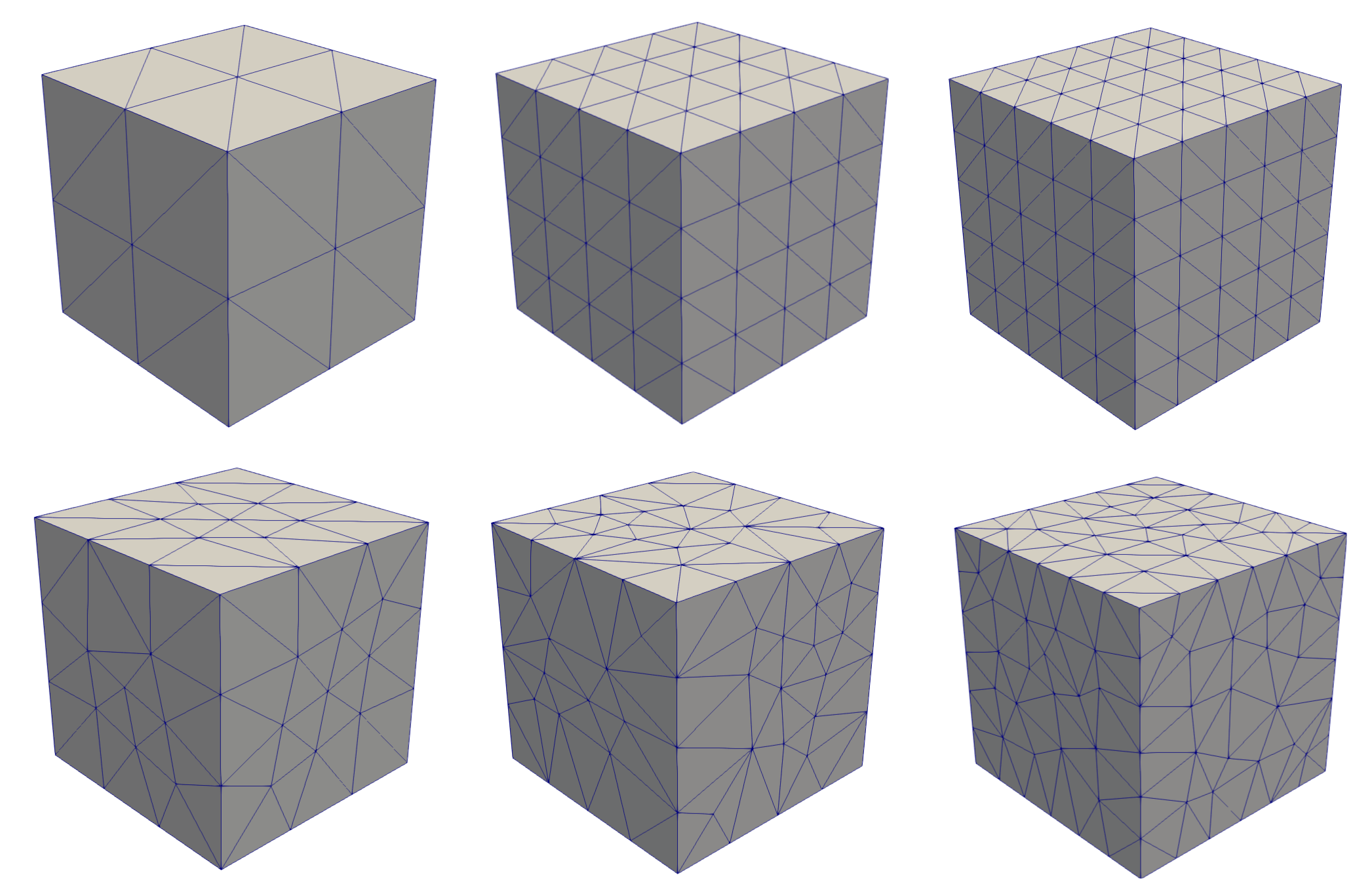}
\end{center}
\vspace*{-0.2in}
\caption{\footnotesize Meshes of the unit cube where the number of elements  $N_{e}$ and the maximum diameter of elements $h$ are given by  $(N_{e},h)= (48,0.866)$, $(384,0.433)$, $(750,0.346)$, for the structured meshes of the first row and $(N_{e},h) = (242, 0.636)$, $(502,0.443)$, $(867,0.363)$, for the unstructured meshes of the second row.} 
\label{CubeStUnstM}
\end{figure}

\begin{table}[!thb]
\tabcolsep=3.0pt
\caption{Convergence rates $r$ and $L^{2}$-errors of the unit cube example: DoF is the number of total degrees of freedom and $(E_{\boldsymbol{U}}, E_{\boldsymbol{K}}, E_{\boldsymbol{P}})=(\|\boldsymbol{U}_{h}-\boldsymbol{U}_{e}\|, \|\boldsymbol{K}_{h}-\boldsymbol{K}_{e}\|, \|\boldsymbol{P}_{h}-\boldsymbol{P}_{e}\|)$ are the $L^{2}$-errors of the approximate solution ($\boldsymbol{U}_{h}$, $\boldsymbol{F}_{h}$, $\boldsymbol{P}_{h}$) with respect to the exact solution ($\boldsymbol{U}_{e}$, $\boldsymbol{F}_{e}$, $\boldsymbol{P}_{e}$) associated to \eqref{cube_Exact}.}
\resizebox{0.8\textwidth}{!}{\begin{minipage}{\textwidth}
\centering
\begin{tabular}{cl|cc|cc|cc||cl|cc|cc|cc}
\toprule
FEM & DoF &$E_{\boldsymbol{U}}$ &  & $E_{\boldsymbol{K}}$ &  & $E_{\boldsymbol{P}}$ &  & FEM & DoF &$E_{\boldsymbol{U}}$ &  & $E_{\boldsymbol{K}}$ &  & $E_{\boldsymbol{P}}$ &  \\
\midrule
\multirow{4}{*}{\rotatebox[origin=c]{90}{\parbox[c]{1.2cm}{\centering $\mathrm{L1N11R1}$}}}
& 735   & 1.82e-2 & \multirow{4}{*}{\rotatebox[origin=c]{90}{\parbox[c]{1.5cm}{\centering $r=2.1$}}}  & 1.21e-1 & \multirow{4}{*}{\rotatebox[origin=c]{90}{\parbox[c]{1.5cm}{\centering $r=1.0$}}}  & 1.67e-1 & \multirow{4}{*}{\rotatebox[origin=c]{90}{\parbox[c]{1.5cm}{\centering $r=0.9$}}} & \multirow{4}{*}{\rotatebox[origin=c]{90}{\parbox[c]{1.2cm}{\centering $\mathrm{L1N11R2}$}}} & 1887   & 1.86e-2 & \multirow{4}{*}{\rotatebox[origin=c]{90}{\parbox[c]{1.5cm}{\centering $r=2.0$}}}  & 1.20e-1 & \multirow{4}{*}{\rotatebox[origin=c]{90}{\parbox[c]{1.5cm}{\centering $r=1.0$}}}  & 1.61e-1 &\multirow{4}{*}{\rotatebox[origin=c]{90}{\parbox[c]{1.5cm}{\centering $r=1.0$}}} \\[2pt]
& 4779   & 4.43e-3 &  & 6.35e-2 &  & 8.97e-2 & & & 13419   & 4.67e-3 &  & 6.23e-2 &  & 8.01e-2 &\\[2pt]
& 8943  & 2.79e-3 &  & 5.12e-2 &  & 7.28e-2 & & & 25593   & 2.96e-3 &  & 5.01e-2 &  & 6.36e-2 &\\[2pt]
\bottomrule
\multirow{4}{*}{\rotatebox[origin=c]{90}{\parbox[c]{1.3cm}{\centering $\mathrm{L1N12R1}$}}}
& 1749   & 1.82e-2 & \multirow{4}{*}{\rotatebox[origin=c]{90}{\parbox[c]{1.5cm}{\centering $r=2.1$}}}  & 1.21e-1 & \multirow{4}{*}{\rotatebox[origin=c]{90}{\parbox[c]{1.5cm}{\centering $r=1.0$}}}  & 1.67e-1 & \multirow{4}{*}{\rotatebox[origin=c]{90}{\parbox[c]{1.5cm}{\centering $r=0.9$}}} & \multirow{4}{*}{\rotatebox[origin=c]{90}{\parbox[c]{1.3cm}{\centering $\mathrm{L1N12R2}$}}} & 2901   & 1.86e-2 & \multirow{4}{*}{\rotatebox[origin=c]{90}{\parbox[c]{1.5cm}{\centering $r=2.0$}}}  & 1.20e-1 & \multirow{4}{*}{\rotatebox[origin=c]{90}{\parbox[c]{1.5cm}{\centering $r=1.0$}}}  & 1.61e-1 &\multirow{4}{*}{\rotatebox[origin=c]{90}{\parbox[c]{1.5cm}{\centering $r=1.0$}}} \\[2pt]
& 11775   & 4.43e-3 &  & 6.35e-2 &  & 8.97e-2 & & & 20415   & 4.67e-3 &  & 6.23e-2 &  & 8.01e-2 &\\[2pt]
& 22188  & 2.79e-3 &  & 5.12e-2 &  & 7.28e-2 & & & 38838   & 2.96e-3 &  & 5.01e-2 &  & 6.36e-2 &\\[2pt]
\bottomrule
\multirow{4}{*}{\rotatebox[origin=c]{90}{\parbox[c]{1.3cm}{\centering $\mathrm{L1N11B1}$}}}
& 1455   & 1.83e-2 & \multirow{4}{*}{\rotatebox[origin=c]{90}{\parbox[c]{1.5cm}{\centering $r=2.1$}}}  & 1.20e-1 & \multirow{4}{*}{\rotatebox[origin=c]{90}{\parbox[c]{1.5cm}{\centering $r=1.0$}}}  & 1.57e-1 & \multirow{4}{*}{\rotatebox[origin=c]{90}{\parbox[c]{1.5cm}{\centering $r=1.1$}}} & \multirow{4}{*}{\rotatebox[origin=c]{90}{\parbox[c]{1.3cm}{\centering $\mathrm{L1N11B2}$}}} & 3399   & 1.86e-2 & \multirow{4}{*}{\rotatebox[origin=c]{90}{\parbox[c]{1.5cm}{\centering $r=2.0$}}}  & 1.20e-1 & \multirow{4}{*}{\rotatebox[origin=c]{90}{\parbox[c]{1.5cm}{\centering $r=1.0$}}}  & 1.63e-1 &\multirow{4}{*}{\rotatebox[origin=c]{90}{\parbox[c]{1.5cm}{\centering $r=1.0$}}} \\[2pt]
& 9963   & 4.51e-3 &  & 6.26e-2 &  & 7.68e-2 & & & 24651   & 4.71e-3 &  & 6.23e-2 &  & 8.19e-2 &\\[2pt]
& 18843  & 2.84e-3 &  & 5.03e-2 &  & 6.05e-2 & & & 47193   & 2.99e-3 &  & 5.00e-2 &  & 6.51e-2 &\\[2pt]
\bottomrule
\multirow{4}{*}{\rotatebox[origin=c]{90}{\parbox[c]{1.3cm}{\centering $\mathrm{L1N12B1}$}}}
& 2469   & 1.83e-2 & \multirow{4}{*}{\rotatebox[origin=c]{90}{\parbox[c]{1.5cm}{\centering $r=2.1$}}}  & 1.20e-1 & \multirow{4}{*}{\rotatebox[origin=c]{90}{\parbox[c]{1.5cm}{\centering $r=1.0$}}}  & 1.57e-1 & \multirow{4}{*}{\rotatebox[origin=c]{90}{\parbox[c]{1.5cm}{\centering $r=1.1$}}} & \multirow{4}{*}{\rotatebox[origin=c]{90}{\parbox[c]{1.3cm}{\centering $\mathrm{L1N12B2}$}}} & 4413   & 1.86e-2 & \multirow{4}{*}{\rotatebox[origin=c]{90}{\parbox[c]{1.5cm}{\centering $r=2.0$}}}  & 2.00e-1 & \multirow{4}{*}{\rotatebox[origin=c]{90}{\parbox[c]{1.5cm}{\centering $r=1.0$}}}  & 1.63e-1 &\multirow{4}{*}{\rotatebox[origin=c]{90}{\parbox[c]{1.5cm}{\centering $r=1.0$}}} \\[2pt]
& 16959   & 4.51e-3 &  & 6.26e-2 &  & 7.68e-2 & & & 31647   & 4.71e-3 &  & 6.23e-2 &  & 8.19e-2 &\\[2pt]
& 32088  & 2.84e-3 &  & 5.03e-2 &  & 6.05e-2 & & & 60438   & 2.99e-3 &  & 5.00e-2 &  & 6.51e-2 &\\[2pt]
\bottomrule
\multirow{4}{*}{\rotatebox[origin=c]{90}{\parbox[c]{1.3cm}{\centering $\mathrm{L1N21R1}$}}}
& 1029   & 1.82e-2 & \multirow{4}{*}{\rotatebox[origin=c]{90}{\parbox[c]{1.5cm}{\centering $r=2.1$}}}  & 1.21e-1 & \multirow{4}{*}{\rotatebox[origin=c]{90}{\parbox[c]{1.5cm}{\centering $r=1.0$}}}  & 1.67e-1 & \multirow{4}{*}{\rotatebox[origin=c]{90}{\parbox[c]{1.5cm}{\centering $r=0.9$}}} & \multirow{4}{*}{\rotatebox[origin=c]{90}{\parbox[c]{1.3cm}{\centering $\mathrm{L1N21R2}$}}} & 2181   & 1.86e-2 & \multirow{4}{*}{\rotatebox[origin=c]{90}{\parbox[c]{1.5cm}{\centering $r=2.0$}}}  & 1.20e-1 & \multirow{4}{*}{\rotatebox[origin=c]{90}{\parbox[c]{1.5cm}{\centering $r=1.0$}}}  & 1.61e-1 &\multirow{4}{*}{\rotatebox[origin=c]{90}{\parbox[c]{1.5cm}{\centering $r=1.0$}}} \\[2pt]
& 6591   & 4.43e-3 &  & 6.35e-2 &  & 8.97e-2 & & & 15231   & 4.67e-3 &  & 6.23e-2 &  & 8.01e-2 &\\[2pt]
& 12288  & 2.79e-3 &  & 5.12e-2 &  & 7.28e-2 & & & 28938   & 2.96e-3 &  & 5.01e-2 &  & 6.36e-2 &\\[2pt]
\bottomrule
\multirow{4}{*}{\rotatebox[origin=c]{90}{\parbox[c]{1.3cm}{\centering $\mathrm{L1N22R1}$}}}
& 2403   & 1.82e-2 & \multirow{4}{*}{\rotatebox[origin=c]{90}{\parbox[c]{1.5cm}{\centering $r=2.1$}}}  & 1.21e-1 & \multirow{4}{*}{\rotatebox[origin=c]{90}{\parbox[c]{1.5cm}{\centering $r=1.0$}}}  & 1.67e-1 & \multirow{4}{*}{\rotatebox[origin=c]{90}{\parbox[c]{1.5cm}{\centering $r=0.9$}}} & \multirow{4}{*}{\rotatebox[origin=c]{90}{\parbox[c]{1.3cm}{\centering $\mathrm{L1N22R2}$}}} & 3555   & 1.86e-2 & \multirow{4}{*}{\rotatebox[origin=c]{90}{\parbox[c]{1.5cm}{\centering $r=2.0$}}}  & 1.20e-1 & \multirow{4}{*}{\rotatebox[origin=c]{90}{\parbox[c]{1.5cm}{\centering $r=1.0$}}}  & 1.61e-1 &\multirow{4}{*}{\rotatebox[origin=c]{90}{\parbox[c]{1.5cm}{\centering $r=1.0$}}} \\[2pt]
& 16179   & 4.43e-3 &  & 6.35e-2 &  & 8.97e-2 & & & 24819   & 4.67e-3 &  & 6.23e-2 &  & 8.01e-2 &\\[2pt]
& 30483  & 2.79e-3 &  & 5.12e-2 &  & 7.28e-2 & & & 47133   & 2.96e-3 &  & 5.01e-2 &  & 6.36e-2 &\\[2pt]
\bottomrule
\multirow{4}{*}{\rotatebox[origin=c]{90}{\parbox[c]{1.3cm}{\centering $\mathrm{L1N21B1}$}}}
& 1749   & 1.83e-2 & \multirow{4}{*}{\rotatebox[origin=c]{90}{\parbox[c]{1.5cm}{\centering $r=2.1$}}}  & 1.20e-1 & \multirow{4}{*}{\rotatebox[origin=c]{90}{\parbox[c]{1.5cm}{\centering $r=1.0$}}}  & 1.57e-1 & \multirow{4}{*}{\rotatebox[origin=c]{90}{\parbox[c]{1.5cm}{\centering $r=1.1$}}} & \multirow{4}{*}{\rotatebox[origin=c]{90}{\parbox[c]{1.3cm}{\centering $\mathrm{L1N21B2}$}}} & 3693   & 1.86e-2 & \multirow{4}{*}{\rotatebox[origin=c]{90}{\parbox[c]{1.5cm}{\centering $r=2.0$}}}  & 1.20e-1 & \multirow{4}{*}{\rotatebox[origin=c]{90}{\parbox[c]{1.5cm}{\centering $r=1.0$}}}  & 1.63e-1 &\multirow{4}{*}{\rotatebox[origin=c]{90}{\parbox[c]{1.5cm}{\centering $r=1.0$}}} \\[2pt]
& 11775   & 4.51e-3 &  & 6.26e-2 &  & 7.68e-2 & & & 26463   & 4.71e-3 &  & 6.23e-2 &  & 8.19e-2 &\\[2pt]
& 22188  & 2.84e-3 &  & 5.03e-2 &  & 6.05e-2 & & & 50538   & 2.99e-3 &  & 5.00e-2 &  & 6.51e-2 &\\[2pt]
\bottomrule
\multirow{4}{*}{\rotatebox[origin=c]{90}{\parbox[c]{1.3cm}{\centering $\mathrm{L1N22B1}$}}}
& 3123   & 1.83e-2 & \multirow{4}{*}{\rotatebox[origin=c]{90}{\parbox[c]{1.5cm}{\centering $r=2.1$}}}  & 1.20e-1 & \multirow{4}{*}{\rotatebox[origin=c]{90}{\parbox[c]{1.5cm}{\centering $r=1.0$}}}  & 1.57e-1 & \multirow{4}{*}{\rotatebox[origin=c]{90}{\parbox[c]{1.5cm}{\centering $r=1.1$}}} & \multirow{4}{*}{\rotatebox[origin=c]{90}{\parbox[c]{1.3cm}{\centering $\mathrm{L1N22B2}$}}} & 5067   & 1.86e-2 & \multirow{4}{*}{\rotatebox[origin=c]{90}{\parbox[c]{1.5cm}{\centering $r=2.0$}}}  & 1.20e-1 & \multirow{4}{*}{\rotatebox[origin=c]{90}{\parbox[c]{1.5cm}{\centering $r=1.0$}}}  & 1.63e-1 &\multirow{4}{*}{\rotatebox[origin=c]{90}{\parbox[c]{1.5cm}{\centering $r=1.0$}}} \\[2pt]
& 21363   & 4.51e-3 &  & 6.26e-2 &  & 7.68e-2 & & & 36051   & 4.71e-3 &  & 6.23e-2 &  & 8.19e-2 &\\[2pt]
& 40383  & 2.84e-3 &  & 5.03e-2 &  & 6.05e-2 & & & 68733   & 2.99e-3 &  & 5.00e-2 &  & 6.51e-2 &\\[2pt]
\bottomrule
\multirow{4}{*}{\rotatebox[origin=c]{90}{\parbox[c]{1.3cm}{\centering $\mathrm{L2N21R2}$}}}
& 2475   & 1.30e-3 & \multirow{4}{*}{\rotatebox[origin=c]{90}{\parbox[c]{1.5cm}{\centering $r=3.0$}}}  & 1.71e-2 & \multirow{4}{*}{\rotatebox[origin=c]{90}{\parbox[c]{1.5cm}{\centering $r=1.9$}}}  & 2.41e-2 & \multirow{4}{*}{\rotatebox[origin=c]{90}{\parbox[c]{1.5cm}{\centering $r=1.9$}}} & \multirow{4}{*}{\rotatebox[origin=c]{90}{\parbox[c]{1.3cm}{\centering $\mathrm{L2N21B2}$}}} & 3987   & 1.30e-3 & \multirow{4}{*}{\rotatebox[origin=c]{90}{\parbox[c]{1.5cm}{\centering $r=3.0$}}}  & 1.71e-2 & \multirow{4}{*}{\rotatebox[origin=c]{90}{\parbox[c]{1.5cm}{\centering $r=2.0$}}}  & 2.40e-2 &\multirow{4}{*}{\rotatebox[origin=c]{90}{\parbox[c]{1.5cm}{\centering $r=2.0$}}} \\[2pt]
& 17043   & 1.66e-4 &  & 4.37e-3 &  & 6.16e-3 & & & 28275   & 1.66e-4 &  & 4.34e-3 &  & 6.13e-3 &\\[2pt]
& 32283  & 8.65e-5 &  & 2.88e-3 &  & 4.33e-3 & & & 53883   & 8.54e-5 &  & 2.79e-3 &  & 3.94e-3 &\\[2pt]
\bottomrule
\multirow{4}{*}{\rotatebox[origin=c]{90}{\parbox[c]{1.3cm}{\centering $\mathrm{L2N12R2}$}}}
& 3195   & 1.30e-3 & \multirow{4}{*}{\rotatebox[origin=c]{90}{\parbox[c]{1.5cm}{\centering $r=3.0$}}}  & 1.71e-2 & \multirow{4}{*}{\rotatebox[origin=c]{90}{\parbox[c]{1.5cm}{\centering $r=1.9$}}}  & 2.41e-2 & \multirow{4}{*}{\rotatebox[origin=c]{90}{\parbox[c]{1.5cm}{\centering $r=1.9$}}} & \multirow{4}{*}{\rotatebox[origin=c]{90}{\parbox[c]{1.3cm}{\centering $\mathrm{L2N12B2}$}}} & 4707   & 1.30e-3 & \multirow{4}{*}{\rotatebox[origin=c]{90}{\parbox[c]{1.5cm}{\centering $r=3.0$}}}  & 1.70e-2 & \multirow{4}{*}{\rotatebox[origin=c]{90}{\parbox[c]{1.5cm}{\centering $r=2.0$}}}  & 2.40e-2 &\multirow{4}{*}{\rotatebox[origin=c]{90}{\parbox[c]{1.5cm}{\centering $r=2.0$}}} \\[2pt]
& 22227   & 1.66e-4 &  & 4.36e-3 &  & 6.16e-3 & & & 33459   & 1.66e-4 &  & 4.34e-3 &  & 6.13e-3 &\\[2pt]
& 42183  & 8.65e-5 &  & 2.88e-3 &  & 4.33e-3 & & & 63783   & 8.54e-5 &  & 2.79e-3 &  & 3.94e-3 &\\[2pt]
\bottomrule
\multirow{4}{*}{\rotatebox[origin=c]{90}{\parbox[c]{1.3cm}{\centering $\mathrm{L2N22R2}$}}}
& 3849   & 1.30e-3 & \multirow{4}{*}{\rotatebox[origin=c]{90}{\parbox[c]{1.5cm}{\centering $r=3.0$}}}  & 1.71e-2 & \multirow{4}{*}{\rotatebox[origin=c]{90}{\parbox[c]{1.5cm}{\centering $r=1.9$}}}  & 2.41e-2 & \multirow{4}{*}{\rotatebox[origin=c]{90}{\parbox[c]{1.5cm}{\centering $r=1.9$}}} & \multirow{4}{*}{\rotatebox[origin=c]{90}{\parbox[c]{1.3cm}{\centering $\mathrm{L2N22B2}$}}} & 5361   & 1.30e-3 & \multirow{4}{*}{\rotatebox[origin=c]{90}{\parbox[c]{1.5cm}{\centering $r=3.0$}}}  & 1.70e-2 & \multirow{4}{*}{\rotatebox[origin=c]{90}{\parbox[c]{1.5cm}{\centering $r=2.0$}}}  & 2.40e-2 &\multirow{4}{*}{\rotatebox[origin=c]{90}{\parbox[c]{1.5cm}{\centering $r=2.0$}}} \\[2pt]
& 26631   & 1.66e-4 &  & 4.36e-3 &  & 6.16e-3 & & & 37863   & 1.66e-4 &  & 4.34e-3 &  & 6.13e-3 &\\[2pt]
& 50478  & 8.65e-5 &  & 2.88e-3 &  & 4.33e-3 & & & 72078   & 8.54e-5 &  & 2.79e-3 &  & 3.94e-3 &\\[2pt]
\bottomrule
\end{tabular} 
\label{CubeConvError}
\end{minipage} }
\end{table}

\begin{figure}[tbh]
\begin{center}
\includegraphics[scale=.5,angle=0]{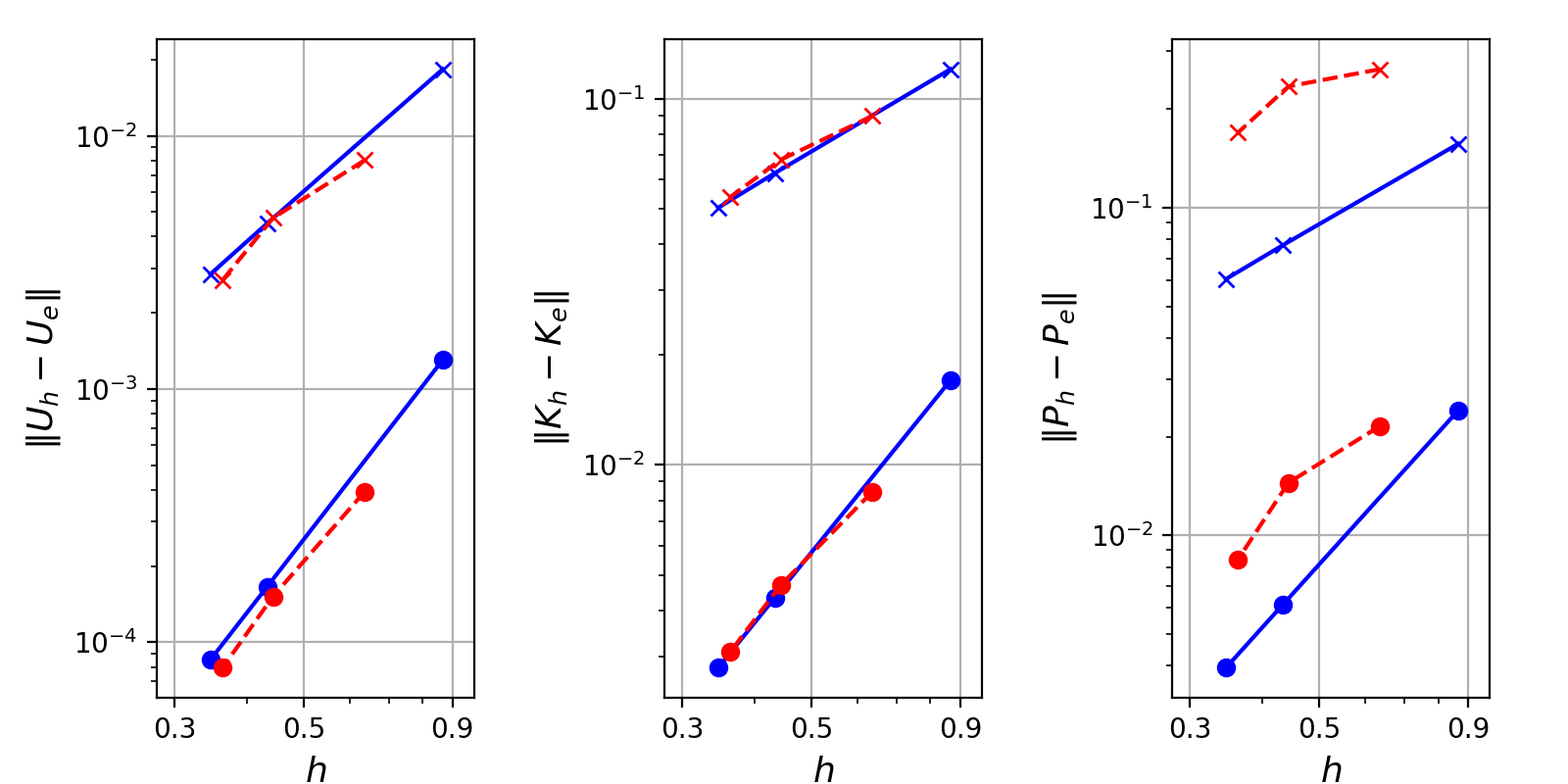}
\end{center}
\vspace*{-0.25in}
\caption{\footnotesize $L^{2}$-errors of displacement $\|\boldsymbol{U}_{h} - \boldsymbol{U}_{e}\|$, displacement gradient $\|\boldsymbol{K}_{h} - \boldsymbol{K}_{e}\|$, and stress $\|\boldsymbol{P}_{h} - \boldsymbol{P}_{e}\|$ associated to the structured meshes (the solid lines) and the unstructured meshes (the dashed lines) of Figure \ref{CubeStUnstM}. The data marked by $\times$ and $\bullet$ are respectively calculated by the first-order elements $\mathrm{L1N21B1}$ and the second-order elements $\mathrm{L2N22B2}$.} 
\label{3DConvStUnst}
\end{figure}

\begin{figure}[tbh]
\begin{center}
\includegraphics[scale=.17,angle=0]{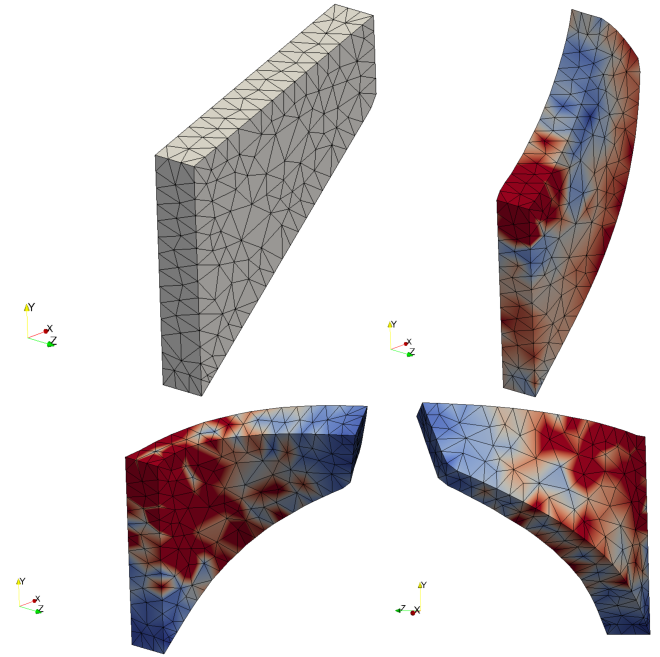}
\end{center}
\vspace*{-0.3in}
\caption{\footnotesize The $3$D Cook-type beam example: The first row shows the geometry (the left panel) and the deformed configuration induced by the uniform in-plane load $F1 = 300\, \mathrm{N/mm^2}$ in the $Y$-direction imposed at the right end (the right panel). The second row shows two different angles of view of the deformed configuration induced by the out-of-plane load $F2 = 600\,\mathrm{N}/\mathrm{mm}^{2}$ in the $Z$-direction applied at the right end of the beam. These results are calculated using the elements $\mathrm{L2N22B2}$ and the underlying mesh has $767$ elements. Colors in these figures depict the distribution of the Frobenius norm of stress.} 
\label{3Dcook_deformed}
\end{figure}

\begin{figure}[tbh]
\begin{center}
\includegraphics[scale=.8,angle=0]{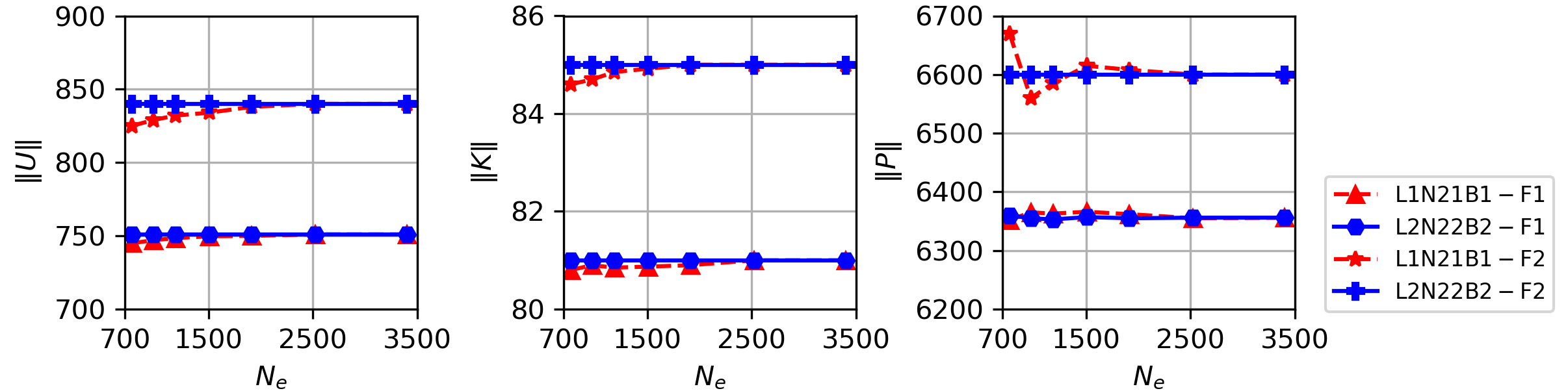}
\end{center}
\vspace*{-0.3in}
\caption{\footnotesize The $L^{2}$-norms of solutions of the $3$D Cook-type beam example associated to the in-plane load $F1 = 300\,\mathrm{N}/\mathrm{mm}^{2}$ and the out-of-plane load $F2 = 600\,\mathrm{N}/\mathrm{mm}^{2}$ versus the number of elements $N_{e}$. The loads $F1$ and $F2$ are imposed at the right end of the beam in the $Y$- and the $Z$-directions, respectively. The elements $\mathrm{L1N21B1}$ and $\mathrm{L2N22B2}$ were used for computing these results.}
\label{3DCOOK_norm}
\end{figure}

%
%

\subsection{Deformation of a Cube}

To study the convergence rates in the $3$D case, we study the $3$D analogue of the plate problem of Section \ref{Ex_2DPlate}. More specifically, we consider the unit cube with the material parameters $\mu=\lambda=1$ and solve the mixed method \eqref{CSFEM23} by using the body force and the boundary conditions that induce the displacement field 
\begin{equation}\label{cube_Exact}
\boldsymbol{U}_{e} = \left[\begin{array}{c} \frac{1}{2}Y^{3} + \frac{1}{2}\sin(\frac{\pi}{2}Y) \\ 0 \\ 0 \end{array} \right].
\end{equation}
Table \ref{CubeConvError} shows $L^{2}$-errors and convergence rates of the solutions of \eqref{CSFEM23}, which are calculated by using different combinations of the $3$D elements of degrees 1 and 2 of Figure \ref{FEO12} and the structured meshes shown in the first row of Figure \ref{CubeStUnstM}. Similar to the $2$D plate example, one observes that the degree of the element for displacement has a significant effect on the overall performance of these mixed finite element methods. Moreover, Table \ref{CubeConvError} suggests that the convergence rates of displacement gradient and stress may not be optimal.

Our numerical results suggest that similar to $2$D cases, $22$ combinations out of $32$ possible combinations of the $3$D elements of Figure \ref{FEO12} are stable. The $10$ unstable cases are the same as those of $2$D cases and are those that do not satisfy the inf-sup conditions \eqref{CSFEM-infsup_Sur} and \eqref{CSFEM-infsup_Sur2}. The extension of the mixed formulation \eqref{CSFEM23} to the $3$D case is straightforward. This is the main advantage of this formulation comparing to the mixed formulation of \citep{AnFSYa2017}.

For the brevity, we consider the choices $\mathrm{L1N21B1}$ and $\mathrm{L2N22B2}$ in the remainder of this work. Figure \ref{3DConvStUnst} depicts the $L^{2}$-errors of approximate solutions corresponding to the structured and unstructured meshes of Figure \ref{CubeStUnstM}. The slopes of curves of the structured meshes are the convergence rates of Table \ref{CubeConvError}. As the $2$D case, these results suggest that mesh irregularities have more impact on the accuracy of approximate stresses.

\subsection{A Near-Incompressible Cook-Type Beam}
Next, we study the $3$D analogue of Cook's membrane under in-plane and out-of-plane loads. The geometry of this problem in the $XY$-plane is similar to that of the $2$D case shown in Figure \ref{CookConfig} with the thickness $10\,\mathrm{mm}$ in the $Z$-direction, see Figure \ref{3Dcook_deformed}. We use the near-incompressible material properties of the $2$D Cook's membrane. 

The configuration in the right panel of the first row of Figure \ref{3Dcook_deformed} is the deformed configuration under a uniform load $F1 = 300\,\mathrm{N}/\mathrm{mm}^{2}$ imposed at the right end of the beam in the $Y$-direction. The second row of Figure \ref{3Dcook_deformed} shows two different angles of view of a deformed configuration due to the out-of-plane load $F2 = 600\,\mathrm{N}/\mathrm{mm}^{2}$ in the $Z$-direction applied at the right end of the beam. These results are computed using $\mathrm{L2N22B2}$ and colors in the deformed configurations depict the distribution of the Frobenius norm of stress. Figure \ref{3DCOOK_norm} shows the convergence of the $L^{2}$-norms of finite element solutions associated to the above in-plane and out-of-plane loads. The elements $\mathrm{L1N21B1}$ and $\mathrm{L2N22B2}$ were used for these computations. Our results suggest that similar to the $2$D case, the $3$D mixed formulation \eqref{CSFEM23} can provide accurate approximations of stress in bending and in the near-compressible regime.

%

\section{Conclusion}\label{Sec_Conc}

We introduced a new mixed formulation for $2$D and $3$D nonlinear elasticity in terms of displacement, displacement gradient, and the first Piola-Kirchhoff stress tensor. We showed that even for hyperelastic solids, this formulation does not correspond to a stationary point of any functional, in general. For obtaining conformal mixed finite element methods based on this formulation, finite element spaces suitable for the $\mathrm{curl}$ and the $\mathrm{div}$ operators are respectively employed for displacement gradient and stress. Discrete displacement gradients and stresses satisfy suitable jump conditions due to these choices.

We studied stability of these mixed finite element methods by writing suitable inf-sup conditions. We examined the performance of these methods for $32$ combinations of $2$D and $3$D simplicial elements of degree $1$ and $2$ and showed that $10$ combinations are not stable as they violate the inf-sup conditions. Several $2$D and $3$D numerical examples were solved to study convergence rates, the effect of mesh distortions, and the performance for bending problems and the near-incompressible regime. These examples suggest that it is possible to achieve the optimal convergence rates and obtain accurate approximations of strains and stresses. Moreover, we did not observe the hourglass instability that may occur in enhanced strain methods.

\appendix

\section{An Abstract Theory for the Galerkin Approximation}\label{Sec_infsup}
In the following, we summarize the general framework for the Galerkin approximation of nonlinear problems introduced in \citep{PoRa1994,CaRa1997}. Let $H:Z\rightarrow Y'$ be a mapping, where $Z$ and $Y$ are Banach spaces with the norms $\|\cdot\|_{Z}$ and $\|\cdot\|_{Y}$, respectively, and $Y'$ is the dual space of $Y$. Also let the linear operator $\mathrm{D}H(u):Z\rightarrow Y'$ be the (Fr\'{e}chet) derivative of $H$ at $u\in Z$, i.e. $\mathrm{D}H(u)z=\frac{d}{ds}|_{s=0} H(u+sz)$, $\forall z\in Z$. The goal is to approximate a regular solution $u\in Z$ of the problem $H(u)=0$, where regular means the derivative of $H$ at $u$ is ``nonzero'' in the sense that the linear mapping $\mathrm{D}H(u)$ is one-to-one and onto. The relation $H(u)=0$ is equivalent to 
\begin{equation}\label{AbsProb}
\langle H(u),y\rangle = 0, \quad \forall y\in Y,
\end{equation} 
where $\langle f,y\rangle:=f(y)$, $\forall f\in Y'$. Given finite element spaces $Z_{h}\subset Z$ and $Y_{h}\subset Y$, a Galerkin approximation of the problem \eqref{AbsProb} reads: Find $u_{h}\in Z_{h}$ such that
\begin{equation}\label{GalAbsProb}
\langle H(u_{h}),y_{h}\rangle = 0, \quad \forall y_{h}\in Y_{h}.
\end{equation} 
To express sufficient conditions for the existence and the convergence of solutions of \eqref{GalAbsProb} as $h\rightarrow 0$, we consider the bilinear form $b:Z\times Y\rightarrow \mathbb{R}$ defined as 
\begin{equation}\label{BiLinF}
b(z,y):=\langle \mathrm{D}H(u)z,y\rangle, \quad \forall z\in Z,~y\in Y.  
\end{equation}
Then, one can show that the following result holds \citep[Theorem 7.1]{CaRa1997}: Roughly speaking, for sufficiently small $h>0$, the problem \eqref{GalAbsProb} has a unique solution $u_{h}$ in a neighborhood of a regular solution $u$ of \eqref{AbsProb} and $u_{h}\rightarrow u$ as $h\rightarrow 0$ if: (i) Any element of $Z$ can be approximated by $Z_{h}$ as $h\rightarrow 0$ (approximibility); (ii) $\dim Z_{h}=\dim Y_{h}$; and (iii) There exists a mesh-independent number $\beta>0$ such that the following inf-sup condition holds:
\begin{equation}\label{infsupCond}
\underset{y_{h}\in Y_{h}}{\inf}\, \underset{z_{h}\in Z_{h}}{\sup} \frac{b(z_{h},y_{h})}{\|z_{h}\|_{Z} \|y_{h}\|_{Y} } \geq \beta>0.
\end{equation}

It is also possible to write a priori and a posteriori estimates for the error $\|u-u_{h}\|_{Z}$ \citep[Theorem 7.1]{CaRa1997}. In particular, the a priori estimate provides an upper bound for $\|u-u_{h}\|_{Z}$ which is proportional to $\beta^{-1}$. If the constant of the inf-sup condition is a mesh-dependent number $\beta_{h}$ such that $\beta_{h}\rightarrow 0$ as $h\rightarrow0$, then $u_{h}$ may converge poorly or diverge as $h\rightarrow0$ even if the inf-sup condition holds for all meshes. Since $u$ is a regular solution of \eqref{AbsProb}, the linearized problem 
\begin{equation*}\label{LinProb}
\langle \mathrm{D} H(u) z,y\rangle = \langle f,y\rangle, \quad \forall y\in Y,
\end{equation*}
has a unique solution $z\in Z$ for any $f\in Y'$. The inf-sup condition \eqref{infsupCond} together with the condition (ii) imply that the discrete linear problem
\begin{equation}\label{DisLinProb}
\langle \mathrm{D} H(u) z_{h},y_{h}\rangle = \langle f,y_{h}\rangle, \quad \forall y_{h}\in Z_{h},
\end{equation}
also has a unique solution $z_{h}\in Z_{h}$ for any $f\in Y'$.

A simple approach to numerically investigate the inf-sup condition \eqref{infsupCond} is as follows: Let $\{\zeta_{i}\}_{i=1}^{n_{Z}}$ and $\{\theta_{i}\}_{i=1}^{n_{Y}}$ respectively be global shape functions for $Z_{h}$ and $Y_{h}$. Then, we have $z_{h}=\sum_{i=1}^{n_{Z}}z_{i}\zeta_{i}$, $\forall z_{h}\in Z_{h}$, and $y_{h}=\sum_{i=1}^{n_{Y}}y_{i}\theta_{i}$, $\forall y_{h}\in Y_{h}$. We associate the vector $\mathbf{z}=(z_{1},\dots,z_{n_{Z}})^{T}\in \mathbb{R}^{n_{Z}}$ ($\mathbf{y}=(y_{1},\dots,y_{n_{Y}})^{T}\in \mathbb{R}^{n_{Y}}$) to $z_{h}$ ($y_{h}$) and define $\|\mathbf{z}\|_{Z}:= \|z_{h}\|_{Z}$ ($\|\mathbf{y}\|_{Y}:= \|y_{h}\|_{Y}$). Assume that there exist symmetric and positive definite matrices $\mathbb{M}^{Z}_{n_{Z}\times n_{Z}}$ and $\mathbb{M}^{Y}_{n_{Y}\times n_{Y}}$ such that 
\begin{align*}
\|\mathbf{z}\|^{2}_{Z}&= (\mathbb{M}^{Z}\mathbf{z})^{T}(\mathbb{M}^{Z}\mathbf{z})= \mathbf{z}^{T}(\mathbb{M}^{Z})^{2}\mathbf{z}, \\
\|\mathbf{y}\|^{2}_{Y}&= (\mathbb{M}^{Y}\mathbf{y})^{T}(\mathbb{M}^{Y}\mathbf{y})= \mathbf{y}^{T}(\mathbb{M}^{Y})^{2}\mathbf{y}.
\end{align*}       
By using the vectors $\mathbf{y}$ and $\mathbf{z}$, the inf-sup condition \eqref{infsupCond} can be expressed in the matrix form
\begin{equation}\label{infsupCondMat}
\underset{\mathbf{y}\in\mathbb{R}^{n_{Y}}}{\inf}\, \underset{\mathbf{z}\in\mathbb{R}^{n_{Z}}}{\sup} \frac{\mathbf{y}^{T}\mathbb{B}\, \mathbf{z}}{\|\mathbf{y}\|_{Y} \|\mathbf{z}\|_{Z}} \geq \beta>0,
\end{equation}
where the matrix $\mathbb{B}_{n_{Y}\times n_{Z}}$ is given by $\mathbb{B}_{ij}=b(\zeta_{j},\theta_{i})$. Recall that the singular values of an arbitrary matrix $\mathbb{M}$ are the square root of the eigenvalues of $\mathbb{M}^{T}\mathbb{M}$. Then, one can show that the inf-sup condition \eqref{infsupCond} holds if and only if the smallest singular value of $\mathbb{M}^{Y}\mathbb{B}\,\mathbb{M}^{Z}$ is positive and bounded from below by a positive constant $\beta$ as $h\rightarrow 0$ \citep[Proposition 3.4.5]{BoBrFo2013}.


\bibliographystyle{unsrtnat}
\bibliography{/home/arzhang/Dropbox/LaTexBibliography/biblio.bib}

\end{document}